\documentclass[10pt]{article}

\RequirePackage{authblk}
\RequirePackage{amsmath}
\RequirePackage{amsfonts}
\RequirePackage{amssymb}
\RequirePackage{enumitem}
\RequirePackage{hyperref}
\RequirePackage{multicol}
\RequirePackage{multirow}
\RequirePackage{algorithm}
\RequirePackage{algpseudocode}
\RequirePackage{graphicx}
\RequirePackage{xcolor}
\RequirePackage{amsthm}
\RequirePackage{amsmath}
\RequirePackage{amsfonts}
\RequirePackage{amssymb}
\RequirePackage{hyperref}
\RequirePackage{algorithm}
\RequirePackage{algpseudocode}
\RequirePackage{graphicx}
\RequirePackage{xcolor}
\RequirePackage{booktabs, multirow}
\RequirePackage{changepage,threeparttable}
\allowdisplaybreaks

\newcommand{\avg}{\mathrm{avg}}
\newcommand{\true}{\mathrm{true}}
\newcommand{\asto}{\xrightarrow{{a.s.}}}
\newcommand{\dto}{\xrightarrow{{d}}}
\newcommand{\pto}{\xrightarrow{{p}}}

\usepackage{dpr_fec,ifthen,dsfont}

\newtheorem{assumption}{Assumption}[section]

\newtheorem{corollary}{Corollary}[section]

\usepackage{isetechreport}
\def\reportyear{23T}
\def\reportno{019}
\def\revisionno{0}
\def\originaldate{\today}
\coraltrue
\cvcrfalse

\begin{document}

\title{Almost-sure convergence of iterates and multipliers in stochastic sequential quadratic optimization}

\author[1]{Frank E.~Curtis\thanks{E-mail: frank.e.curtis@lehigh.edu}}
\author[1]{Xin Jiang\thanks{E-mail: xjiang@lehigh.edu}}
\author[1]{Qi Wang\thanks{E-mail: qiw420@lehigh.edu}}
\affil[1]{Department of Industrial and Systems Engineering, Lehigh University}
\titlepage

\maketitle

\begin{abstract}
      Stochastic sequential quadratic optimization (SQP) methods for solving continuous optimization problems with nonlinear equality constraints have attracted attention recently, such as for solving large-scale data-fitting problems subject to nonconvex constraints.  However, for a recently proposed subclass of such methods that is built on the popular stochastic-gradient methodology from the unconstrained setting, convergence guarantees have been limited to the asymptotic convergence of the expected value of a stationarity measure to zero.  This is in contrast to the unconstrained setting in which almost-sure convergence guarantees (of the gradient of the objective to zero) can be proved for stochastic-gradient-based methods.  In this paper, new almost-sure convergence guarantees for the primal iterates, Lagrange multipliers, and stationarity measures generated by a stochastic SQP algorithm in this subclass of methods are proved.  It is shown that the error in the Lagrange multipliers can be bounded by the distance of the primal iterate to a primal stationary point plus the error in the latest stochastic gradient estimate.  It is further shown that, subject to certain assumptions, this latter error can be made to vanish by employing a running average of the Lagrange multipliers that are computed during the run of the algorithm.  The results of numerical experiments are provided to demonstrate the proved theoretical guarantees.
\end{abstract}

\section{Introduction}\label{sec.introduction}

In this paper, we study convergence guarantees that can be offered for a stochastic algorithm for solving continuous optimization problems with nonlinear equality constraints. Such problems arise in a variety of important areas throughout science and engineering, including optimal control, PDE-constrained optimization, network optimization, and resource allocation \cite{Ber98,Bet10,KS92,RDW10}, and have recently arisen in new and interesting modern application areas such as constrained deep neural network training (e.g., physics-informed learning \cite{CDG+22,KKL+21} where one can impose hard constraints rather than merely define the loss function to minimize residual errors \cite{BCRZ21}). In certain instances of such problems, the objective function can be defined as an expectation over a random variable argument.  In this context, the useful features of the algorithm that we study are that it is applicable when only (unbiased) stochastic estimates of the gradient of the objective function are tractable to obtain during the optimization, while at the same time the algorithm can exploit exact constraint function and derivative values that are tractable to obtain.

Sequential quadratic optimization (commonly known as SQP) methods are a popular and powerful class of derivative-based algorithms for solving continuous constrained optimization problems.  SQP methods that are stochastic in nature (due to their use of stochastic gradient estimates in place of true gradients of the objective) have been investigated recently for solving problems of the aforementioned type, namely, ones for which the objective is defined by an expectation.  For one subclass of stochastic SQP methods developed in \cite{BCRZ21,BCOR23,COR23,CRZ23} that is based on the classical stochastic-gradient methodology (known more broadly as stochastic approximation) \cite{RM51}, convergence guarantees have been limited to asymptotic convergence of the expected value of a stationarity measure to zero; see, e.g., \cite[Corollary~3.14]{BCRZ21}.  In this paper, we present new analyses of situations in which, for an instance of a method in this subclass, almost-sure convergence guarantees of the primal iterates, Lagrange multipliers, and a stationarity measure can be guaranteed.  (For a straightforward analysis, we focus on a simplified variant of the algorithm from \cite{BCRZ21}.)  This brings the analysis of this subclass of methods in line with analyses that show almost-sure convergence of stochastic-gradient-based methods for the unconstrained setting; see \S\ref{sec.background}.  Other stochastic SQP-based methods for solving continuous optimization problems with nonlinear equality constraints have been proposed recently (see, e.g., \cite{BBZ23,BSYZ23,BXZ23,FNMK22,NAK23a,NAK23,OBN21,QK23}), and for these methods different types of convergence guarantees have been proved due to the fact that they impose different (often stronger) requirements on the stochastic gradient estimates.  Therefore, our work in this paper is distinct from analyses for these other methods, although we contend that the tools employed in this paper might be applicable when analyzing other related algorithms as well.  One article that presents results comparable to ours is \cite{NM22}, although in that article the algorithm and assumptions are different than our setting.

Our first contribution is an analysis of situations in which almost-sure convergence of the primal iterates of a stochastic SQP method can be guaranteed.  Our main assumption for this analysis can be viewed as a generalization of the Polyak-Łojasiewicz (PL) condition used in unconstrained optimization.

Our second contribution is an analysis of convergence of the sequence of Lagrange multipliers generated by a stochastic SQP method.  The convergence behavior of the Lagrange multipliers is of interest for several reasons. For one thing, Lagrange multipliers are used for common certificates of stationarity. They can also play an essential role in sensitivity analysis and active-set identification when solving inequality-constrained problems. Since various algorithms for solving inequality-constrained problems employ algorithms for solving equality-constrained subproblems, it is important to provide convergence guarantees for the Lagrange multipliers when solving equality-constrained problems.  In our analysis, we first show that the expected error in the Lagrange multiplier computed in any iteration of our stochastic SQP method of interest can be bounded by the distance of the primal iterate to a primal stationary point plus a term related to the error in the latest stochastic gradient estimate.  Such a result is natural, and as in other settings of statistical estimation it suggests that better Lagrange multiplier estimates can be obtained through averaging the multipliers obtained in each iteration.  We consider such an approach as well; in particular, we show conditions under which averaging can cause the expected error to vanish asymptotically.  (Our approach to averaging is related to the ergodic convergence analysis of optimization algorithms; see, e.g., \cite{PJ92,CP16,JV23}. It is worthwhile to emphasize that it is different from dual averaging ideas that have been developed for acceleration \cite{Nes05}.)

Under our combined assumptions, the stochastic SQP method that we analyze possesses almost-sure convergence guarantees of the primal iterates, Lagrange multipliers, and a stationarity measure while only employing stochastic gradient estimates.  To illustrate the practical relevance of our theoretical contributions, we provide the results of numerical experiments. For example, we demonstrate situations when averaging of the Lagrange multipliers results in more accurate stationarity measures, which, as mentioned, is useful in practice for recognizing (approximate) stationarity and other reasons.

\subsection{Notation}\label{sec.notation}

We use $\R{}$ to denote the set of real numbers and $\R{}_{\geq r}$ (resp., $\R{}_{>r}$) to denote the set of real numbers greater than or equal to (resp., greater than) $r \in \R{}$. We use $\R{n}$ to denote the set of $n$-dimensional real vectors and $\R{m \times n}$ to denote the set of $m$-by-$n$-dimensional real matrices.  We use $\mathbb{S}^n$ to denote the set of symmetric matrices in $\R{n \times n}$, $\mathbb{S}_{\succeq0}^n$ to denote the set of positive semidefinite matrices in~$\mathbb{S}^n$, and $\mathbb{S}_{\succ0}^n$ to denote the set of positive definite matrices in $\mathbb{S}_{\succeq0}^n$.  We use $\N{}$ to denote the positive integers, and for any $k \in \N{}$ define $[k] := \{1,\dots,k\}$.

Given a real matrix $A \in \R{n \times m}$, we use $\sigma_{\min}(A)$ to denote the smallest singular value of $A$ and $\sigma_{\max}(A) = \|A\|_2$ to denote the largest singular value of $A$, i.e., the spectral norm of $A$. Given such a matrix $A$ with full column rank, we denote by $A^\dagger = (A^TA)^{-1}A^T$ the Moore--Penrose pseudoinverse of $A$.  For future reference, we state the following lemma pertaining to pseudoinverses.

\blemma[see\text{\cite[Theorem~4.1]{Wed73}}]\label{lem.pinv}
  If $A \in \R{n \times m}$ and $B \in \R{n \times m}$ have full column rank, then $\|A^\dagger-B^\dagger\|_2 \leq \|A^\dagger\|_2 \|B^\dagger\|_2 \|A-B\|_2$.
\elemma

For a real-valued sequence $\{x_k\}$ (of numbers, vectors, or matrices), we use $\{x_k\} \subset \R{n}$ to indicate that $x_k \in \R{n}$ for all $k \in \N{}$.  Similarly, for a sequence of random variables $\{X_k\}$, we use $\{X_k\} \subset \R{n}$ to indicate that $X_k \in \R{n}$ for all $k \in \N{}$, which in turn means that, for all $k \in \N{}$, a realization/outcome of~$X_k$ is an element of $\R{n}$.  Generally, we use a capital letter to denote a random variable and the corresponding lower-case letter to denote a realization of the random variable.  For example, a stochastic objective gradient estimator at iteration $k \in \N{}$ is denoted as $G_k$, of which a realization is written as $g_k$.

Let $\{V_k\} \subset \R{n \times m}$ be a sequence of random variables and $V$ be a random variable all defined with respect to a probability space $(\Omega, \Fcal, \P)$; i.e., an realization $\omega \in \Omega$ defines a realization $V_k(\omega)$ for any $k \in \N{}$ and a realization $V(\omega)$.  The sequence~$\{V_k\}$ converges in distribution to $V$ if and only if the cumulative distribution functions (CDFs) of the elements of $\{V_k\}$ converge pointwise to the CDF of~$V$ as $k \to \infty$; e.g., for such $\{V_k\} \subset \R{n}$ and a matrix $\Sigma \in \mathbb{S}^n$, we write $\{V_k\} \dto \Ncal(0, \Sigma)$ to indicate that $\{V_k\}$ converges in distribution to a multivariate normal random vector with mean zero and covariance matrix $\Sigma$.  The sequence $\{V_k\}$ converges in probability to $V$ as $k \to \infty$, which we denote by $\{V_k\} \pto V$, if and only if for any $\epsilon \in \R{}_{>0}$ one finds
\bequationNN
  \lim_{k \to \infty} \P[\|V_k - V\|_2 > \epsilon] = 0.
\eequationNN
Finally, $\{V_k\}$ converges almost-surely to $V$ as $k \to \infty$, which we denote by $\{V_k\} \asto V$, if and only if there exists $\Omega_0 \subset \Omega$ with $\P(\Omega_0) = 0$ such that
\bequationNN
  \lim_{k \to \infty} V_k(\omega) = V(\omega)\ \ \text{for all}\ \ \omega \in \Omega \setminus \Omega_0.
\eequationNN

We use $\mathbf{1}_A$ to denote an indicator random variable for the event $A$, which takes the value 1 if event $A$ occurs and takes the value 0 otherwise.

\subsection{Mathematical Background}\label{sec.background}

In the classical article by Robbins and Monro \cite{RM51}, it is shown that a straightforward approach of stochastic approximation for solving an equation with a unique root $x_\star$ can lead to convergence in probability of the iterate sequence.  Specifically, under certain basic assumptions, as long as the prescribed sequence of \emph{step sizes} $\{\alpha_k\}$ employed by the algorithm is unsummable, but square summable (e.g., $\alpha_k = 1/k$ for all $k \in \N{}$), it can be shown that the generated sequence of solution estimates $\{X_k\}$ satisfies
\bequation\label{eq.robbins}
  \lim_{k \to \infty} \E[(X_k - x_\star)^2] = 0,
\eequation
which in turn implies that $\{X_k\} \pto x_\star$.  Cast into the context of minimizing a smooth, potentially nonconvex objective $f : \R{n} \to \R{}$ with a stochastic-gradient method, these same principles can be used to prove (see, e.g., \cite{BCN18})
\bequationNN
  \lim_{k \to \infty} \E[\|\nabla f(X_k)\|_2^2] = 0.
\eequationNN

In a later article by Robbins and Siegmund \cite{RS71}, certain conditions are shown to offer an almost-sure convergence guarantee for a sequence generated by stochastic approximation. In particular, the article first proves a general theorem, which we state in a slightly simplified form for our purposes.

\blemma[see\text{\cite[Theorem~1]{RS71}}]\label{lem.rs}
  Let $(\Omega,\Fcal,\P)$ be a probability space and let $\{\Fcal_k\}$ with $\Fcal_k \subseteq \Fcal_{k+1}$ for all $k \in \N{}$ be a sequence of sub-$\sigma$-algebras of $\Fcal$.  Let $\{R_k\}$, $\{P_k\}$, and $\{Q_k\}$ be sequences of nonnegative random variables such that, for all $k \in \N{}$, the random variables $R_k$, $P_k$, and $Q_k$ are $\Fcal_k$-measurable.  If $\sum_{k=1}^\infty Q_k < \infty$ and, for all $k \in \N{}$, one has
  \bequationNN
    \E[R_{k+1} | \Fcal_k] \leq R_k - P_k + Q_k,
  \eequationNN
  then, almost-surely, $\sum_{k=1}^\infty P_k < \infty$ and $\displaystyle \lim_{k \to \infty} R_k$ exists and is finite.
\elemma

\noindent
Applied to the context of stochastic approximation for solving an equation, it is shown in \cite{RS71} that under certain assumptions (which we omit for brevity) the conclusion in \eqref{eq.robbins} can be strengthened for the same algorithm to
\bequationNN
  \P\left[ \lim_{k \to \infty} X_k = x_\star \right] = 1,
\eequationNN
which is to say that $\{X_k\} \asto x_\star$.  One can derive a similar such result in the context of minimizing a smooth, potentially nonconvex objective function $f : \R{n} \to \R{}$ with a stochastic-gradient method.  For example, under a related set of conditions, it has been shown by Bertsekas and Tsitsiklis \cite{BT00} that a stochastic-gradient method can be guaranteed to yield
\bequation\label{eq.bertsekas}
  \P\left[ \lim_{k \to \infty} \nabla f(X_k) = 0 \right] = 1,\ \ \text{i.e.,}\ \ \{\nabla f(X_k)\} \asto 0.
\eequation

Our main contributions in this paper are also almost-sure convergence guarantees, but for a stochastic SQP method in the context of nonlinear-equality-constrained optimization.  In particular, we show almost-sure convergence guarantees for the primal iterates, Lagrange multipliers, and stationarity measures for a simplified variant of the algorithm from \cite{BCRZ21}.  Some of our almost-sure convergence guarantees for the Lagrange multipliers computed by our method of interest pertains to an averaged sequence.  For our analysis of this sequence, we require two key results that are known from the literature.  The first result that we need is the central limit theorem (CLT) for a multidimensional martingale difference triangular array stated as Lemma~\ref{lem.mart_clt}.

\blemma[Multidimensional martingale central limit theorem]\label{lem.mart_clt}
  Let $\{(\xi_{k,i}, \Fcal_{k,i})\}_{k \in \N{},i \in [k]}$ be an $n$-dimensional martingale difference triangular array, i.e., with an initial generating $\sigma$-algebra $\Fcal_{k,1}$ for all $k \in \N{}$, one has
  \benumerate
    \item[(i)] $\Fcal_{k,i} = \sigma(\xi_{k,1},\dots,\xi_{k,i-1})$ for all $k \in \N{}$ and $i \in \{2,\dots,k\}$ and
    \item[(ii)] $\E[\xi_{k,i} | \Fcal_{k,i}] = 0$ for all $k \in \N{}$ and $i \in [k]$.
  \eenumerate
  If the array has the properties that
  \begin{subequations} \label{eq.mart_clt_cond}
    \begin{align}
      &\text{$\xi_{k,i}$ is square-integrable, i.e., $\E[\|\xi_{k,i}\|_2^2] < \infty$ for all $(k,i) \in \N{} \times [k]$,} \label{eq.mart_clt_int} \\
      &\text{$\left\{\sum_{i=1}^k \E[\|\xi_{k,i}\|_2^2 \mathbf{1}_{\{\|\xi_{k,i}\|_2 > \delta\}} | \Fcal_{k,i}]\right\} \pto 0$ for all $\delta \in \R{}_{>0}$, and} \label{eq.mart_clt_lindeberg} \\
      &\text{$\left\{\sum_{i=1}^k \E[\xi_{k,i} \xi_{k,i}^T | \Fcal_{k,i}]\right\} \pto \Sigma$ for some $\Sigma \in \mathbb{S}^n$}, \label{eq.mart_clt_var}
    \end{align}
  \end{subequations}
  then $\displaystyle \left\{\sum_{i=1}^k \xi_{k,i}\right\} \dto \Ncal(0, \Sigma)$.
\elemma

\noindent
The multidimensional martingale CLT in Lemma~\ref{lem.mart_clt} can be derived by applying the one-dimensional martingale CLT (see, e.g., \cite[Theorem~2.3]{McL74} and \cite[Corollary~3.1]{HH14}) to $\xi_{k,i}^a := a^T \xi_{k,i}$ for arbitrary $a \in \R{n}$ since $\{\sum_{i=1}^k \xi_{k,i}^a\} \dto \Ncal(0, a^T\Sigma a)$ for any given $a \in \R{n}$ implies that $\{\sum_{i=1}^k \xi_{k,i}\} \dto \Ncal(0,\Sigma)$; see, e.g.,~\cite[Exercise~3.10.8]{Dur19}.  A similar result is used in \cite{LXY23}.  We refer to \eqref{eq.mart_clt_lindeberg} as Lindeberg's condition, as is common in the literature.

The second key result that we need is the following, which we refer to as a moment convergence result.  It follows, e.g., by \cite[Theorem 4.5.2]{Chu03}.

\blemma[Moment convergence]\label{lem.moment_convergence}
  Let $\{X_k\} \subset \R{n}$ be a sequence of random vectors such that $\{X_k\} \dto X$ and $\displaystyle \sup_{k\in\N{}} \E[\|X_k\|_2^\Theta] < \infty$ for some $\Theta \in \R{}_{>0}$.  Then, for all $\theta \in (0,\Theta)$, one finds that $\displaystyle \lim_{k\to\infty} \E[\|X_k\|_2^\theta] = \E[\|X\|_2^\theta]$.
\elemma

\subsection{Outline}\label{sec.outline}

In Section~\ref{sec.algorithm}, we present formally the continuous optimization problem of interest and the stochastic SQP algorithm for solving it that we analyze in the remainder of the paper.  We also provide preliminary assumptions that we make about the problem and the algorithm, and state basic properties of the algorithm that transfer from the prior work in \cite{BCRZ21}. In Section~\ref{sec.primal_iterates}, we prove convergence results for the primal iterates generated by the algorithm. In Section~\ref{sec.multipliers}, we prove convergence results for Lagrange multiplier sequences that are generated by the algorithm. The results of numerical experiments are provided in Section~\ref{sec.numerical} and concluding remarks are offered in Section~\ref{sec.conclusion}.

\section{Problem and Algorithm Descriptions}\label{sec.algorithm}

The algorithm that we study is designed to solve
\bequation\label{prob.opt}
  \min_{x \in \R{n}}\ f(x)\ \text{subject to (s.t.)}\ c(x) = 0\ \text{with}\ f(x) = \E_\iota[F(x,\iota)],
\eequation
where $f : \R{n} \to \R{}$ and $c : \R{n} \to \R{m}$ are continuously differentiable, $\iota$ is a random variable with associated probability space $(\Omega_\iota, \Fcal_\iota, \P_\iota)$, $F : \R{n} \times \Omega_\iota \to \R{}$, and $\E_\iota [\cdot]$ denotes expectation taken with respect to $\P_\iota$.  Given an initial point $x_1 \in \R{n}$, any run of the algorithm generates a sequence of iterates $\{x_k\} \subset \R{n}$, i.e., a realization of a stochastic process $\{X_k\} \subset \R{n}$.  We make the following assumption about problem~\eqref{prob.opt} and the generated sequence of iterates.

\begin{assumption}\label{ass.opt}
  There exists open convex $\Xcal \subseteq \R{n}$ containing $\{X_k\} \subset \R{n}$ such that the following hold.  The objective function $f : \R{n} \to \R{}$ is continuously differentiable and bounded below over~$\Xcal$ and its gradient function $\nabla f : \R{n} \to \R{n}$ is bounded and Lipschitz continuous over~$\Xcal$.  The constraint function $c : \R{n} \to \R{m}$ $($where $m \leq n$$)$ is Lipschitz continuous, continuously differentiable, and bounded over~$\Xcal$, its Jacobian function $\nabla c^T : \R{n} \to \R{m \times n}$ is Lipschitz continuous over $\Xcal$, and $\nabla c(x)^T$ has full row rank for all $x \in \cal X$ with singular values that are bounded below uniformly by a positive constant over~$\Xcal$.
\end{assumption}

\noindent
It would be possible to loosen Assumption~\ref{ass.opt} to require only that~$\Xcal$ contains $\{X_k\}$ almost surely.  However, since this would require repeated references to probability-one events throughout our analysis without strengthening our conclusions substantially, we do not bother with this level of generality.

Under Assumption~\ref{ass.opt}, there exists a tuple of constants, which we denote by $(\kappa_\Xcal, f_{\inf}, \kappa_{\nabla f}, \kappa_c, \kappa_{\nabla c}, r) \in \R{}_{>0} \times \R{} \times \R{}_{>0} \times \R{}_{>0} \times \R{}_{>0} \times \R{}_{>0}$, such that
\bequation\label{eq.ass_opt}
  \begin{gathered}
    \|x\|_2 \leq \kappa_\Xcal,\ f(x) \geq f_{\inf},\ \|\nabla f(x)\|_2 \leq \kappa_{\nabla f},\\ \|c(x)\|_2 \leq \kappa_c,\ \|\nabla c(x)\|_2 \leq \kappa_{\nabla c},\ \text{and}\ \sigma_{\min}(\nabla c(x)) \geq r
  \end{gathered}
\eequation
for all $x \in \cal X$, and there exists $(L_{\nabla f}, L_c, \Gamma) \in \R{}_{>0} \times \R{}_{>0} \times \R{}_{>0}$ such that
\bequation\label{eq.Lipschitz}
  \begin{gathered}
    \|\nabla f(x) - \nabla f(\xbar)\|_2 \leq L_{\nabla f} \|x - \xbar\|_2,\ \|c(x) - c(\xbar)\|_2 \leq L_c \|x - \xbar\|_2, \\
    \text{and}\ \ \|\nabla c(x) - \nabla c(\xbar)\|_2 \leq \Gamma \|x - \xbar\|_2
  \end{gathered}
\eequation
for all $(x,\xbar) \in \Xcal \times \Xcal$.  The Lagrangian function $\Lcal : \R{n} \times \mathbb R^m \to \R{}$ corresponding to problem~\eqref{prob.opt} is defined by $\Lcal(x,y) = f(x) + c(x)^Ty$, and the first-order stationarity conditions for~\eqref{prob.opt} are given by
\bequation\label{eq.kkt}
  0 = \bbmatrix \nabla_x \Lcal(x,y) \\ \nabla_y \Lcal(x,y) \ebmatrix = \bbmatrix \nabla f(x) + \nabla c(x)y \\ c(x) \ebmatrix.
\eequation
We refer to any $(x,y) \in \R{n} \times \R{m}$ satisfying \eqref{eq.kkt} as a stationary point for problem~\eqref{prob.opt}.  In addition, we refer to any $x \in \R{n}$ such that there exists $y \in \R{m}$ with $(x,y)$ satisfying \eqref{eq.kkt} as a primal stationary point for \eqref{prob.opt}.

Let us now describe the algorithm whose convergence properties are the subject of our study.  We state the algorithm in terms of a realization of the quantities it generates.  At an iterate $x_k \in \R{n}$, the algorithm generates a stochastic gradient estimate $g_k \approx \nabla f(x_k) \in \R{n}$, and makes use of $h_k \in \mathbb{S}^n$ (see upcoming Assumptions~\ref{ass.g} and~\ref{ass.H} about these quantities). Given $g_k$ and~$h_k$, a direction $d_k \in \R{n}$ is computed by solving
\bequation\label{prob.d}
  \min_{d \in \R{n}}\ g_k^Td + \tfrac12 d^Th_kd\ \st\ c(x_k) + \nabla c(x_k)^Td = 0.
\eequation
Under upcoming Assumption~\ref{ass.H} (which includes that $d^Th_kd > 0$ for all $d \in \Null(\nabla c(x_k)^T)$), the solution $d_k$ of \eqref{prob.d} as well as a Lagrange multiplier for the constraints, call it $y_k \in \R{m}$, can be obtained by solving
\bequation\label{eq.system}
  \bbmatrix h_k & \nabla c(x_k) \\ \nabla c(x_k)^T & 0 \ebmatrix \bbmatrix d_k \\ y_k \ebmatrix = - \bbmatrix g_k \\ c_k \ebmatrix.
\eequation

Upon computation of the search direction in iteration $k \in \N{}$, the algorithm selects a step size $\alpha_k \in (0,1]$.  Specifically, with a \emph{merit parameter} $\tau \in \R{}_{>0}$, \emph{ratio parameter} $\xi \in \R{}_{>0}$ (see upcoming Assumption \ref{ass.tau_xi}), and $(L_{\nabla f},\Gamma)$ from~\eqref{eq.Lipschitz}, the algorithm that we analyze selects the step size for all $k \in \N{}$ as
\bequation\label{eq.alpha}
  \alpha_k \gets \tfrac{\beta_k \tau \xi}{\tau L_{\nabla f} + \Gamma}\ \ \text{for some} \ \beta_k \in (0,1],
\eequation
where $\{\beta_k\}$ is unsummable, but square-summable. This choice of step sizes means that the algorithm that we analyze is a simplified variant of the algorithm from \cite{BCRZ21}; see \cite{BCRZ21} and below for further discussion about the merit parameter, ratio parameter, and the particular formula for the step size stated in \eqref{eq.alpha}.  (We conjecture that our ultimate conclusions also hold for the original variant of the algorithm from \cite{BCRZ21} in the event considered in Section~3.2.1 of that paper, but for our aims in this paper we merely consider the simplified variant presented here in order to allow our analysis not to be obscured by auxiliary details.) One aspect that makes our variant a simplified one is that we assume that the merit and ratio parameters are initialized to values that are sufficiently small such that the algorithm does not need to update them adaptively.  This means that, for our analysis, we can consider the merit function (with \emph{fixed} merit parameter $\tau$) $\phi_\tau : \R{n} \to \R{}$ defined by $\phi_\tau(x) = \tau f(x) + \|c(x)\|_1$.  (Any convex norm can be used for the constraint violation.  We use the $\ell_1$-norm for consistency with \cite{BCRZ21}.)  A local approximation of $\phi_\tau$ at $x \in \R{n}$ can be defined through a function $q_\tau \colon \R{n} \times \R{n} \times \mathbb{S}^n \times \R{n} \to \R{}$; specifically, with $g \in \R{n}$ and $h \in \mathbb{S}^n$, a local approximation of $\phi_\tau$ at $x$ as a function of $d$ is given by
\bequation\label{eq.merit_model}
  q_\tau(x,g,h,d) = \tau (f(x) + g^Td + \tfrac12 \max\{d^Thd,0\}) + \|c(x) + \nabla c(x)^T d\|_1.
\eequation
The reduction in $q_\tau(x,g,h,\cdot)$ offered by $d$ with $c(x) + \nabla c(x)^Td = 0$ can be defined through $\Delta q_\tau \colon \R{n} \times \R{n} \times \R{n \times n} \times \R{n} \to \R{}$ defined by
\begin{align}
  \Delta q_\tau(x,g,h,d) &= q_\tau(x,g,h,0) - q_\tau(x,g,h,d) \nonumber \\
  &= -\tau (g^Td + \tfrac12 \max\{d^Thd,0\}) + \|c(x)\|_1.\label{eq.merit_model_reduction}
\end{align}
Consistent with \cite[Section~3.2.1]{BCRZ21}, we make the following assumption.

\begin{assumption}\label{ass.tau_xi}
  The parameters $\xi \in \R{}_{>0}$ and $\tau \in \R{}_{>0}$ are chosen such that, for some $\nu \in (0,1)$ and all $k \in \N{}$ in any run, the following hold.
  \benumerate
    \item[(i)] $\xi \leq \xi_k^{\rm trial}$, where one defines
    \bequationNN
      \xi_k^{\rm trial} := \bcases \infty & \text{if $d_k = 0$} \\ \tfrac{\Delta q_\tau(x_k, g_k, h_k, d_k)}{\tau \|d_k\|_2^2} & \text{otherwise;} \ecases
    \eequationNN
\item[(ii)] $\tau \leq \tau_k^{\rm trial, true}$, where, with $d_k^\true$ being the solution of \eqref{prob.d} if $g_k$ were replaced by $\nabla f(x_k)$ and $\rho_k := \nabla f(x_k)^Td_k^\true + \max\{(d_k^\true)^Th_kd_k^\true, 0\}$, one defines
    \begin{equation*}
      \tau_k^{\rm trial,\, true} := \begin{cases} \infty & \text{if $\rho_k \leq 0$} \\ \displaystyle \tfrac{(1 - \nu)\|c(x_k)\|_1}{\rho_k} & \text{otherwise.} \end{cases}
    \end{equation*}
  \eenumerate
\end{assumption}

\noindent
The definition of $\tau_k^{\rm trial,\, true}$ and Assumption~\ref{ass.tau_xi} ensure for all $k \in \N{}$ that
\bequation \label{eq.merit_reduction_lbnd}
  \Delta q_\tau(x_k, \nabla f(x_k), h_k, d_k^{\true}) \ge \thalf \tau
\max \{(d_k^{\true})^T h_k d_k^{\true}, 0\} + \nu \|c_k\|_1.
\eequation
Overall, the step-size choice in \eqref{eq.alpha} along with Assumption~\ref{ass.tau_xi} results in a simplified variant of the adaptive strategy presented in \cite{BCRZ21}, and one finds from~\cite{BCRZ21} that this setting ensures that the step from $x_k$ to $x_{k+1}$ yields a sufficient-decrease-type property that is relevant for our analysis; see Lemma~\ref{lem.sufficient_decrease} below.

The stochastic SQP method that we study is stated in Algorithm \ref{alg.sqp}.

\begin{algorithm}[ht]
  \caption{Stochastic SQP}
  \label{alg.sqp}
  \begin{algorithmic}[1]
    \Require $x_1 \in \R{n}$; $(L_{\nabla
    f},\Gamma) \in \R{}_{>0} \times \R{}_{>0}$ (see~\eqref{eq.Lipschitz}); $(\tau,\xi) \in \R{}_{>0} \times \R{}_{>0}$ satisfying Assumption~\ref{ass.tau_xi}; $\{\beta_k\} \subset (0,1]$ that is unsummable, but square-summable.
    \For{$k\in\N{}$}
      \State Compute $g_k$ (see Assumption~\ref{ass.g}). \label{step.g}
      \State Compute $h_k$ (see Assumption~\ref{ass.H}).
      \State Compute $(d_k,y_k)$ by solving \eqref{prob.d}.
      \State Set $\alpha_k$ by \eqref{eq.alpha}.
      \State Set $x_{k+1}\gets x_k + \alpha_k d_k$.
    \EndFor
  \end{algorithmic}
\end{algorithm}

All that remains before presenting our analysis are to articulate some final assumptions and state a few key relations from \cite{BCRZ21}.  Considering the outcomes of entire runs of Algorithm~\ref{alg.sqp}, henceforth we consider the probability space $(\Omega, \Fcal, \P)$, where $\Omega := \prod_{k=1}^\infty \Omega_\iota$.  In this manner, each realization of a run of the algorithm can be associated with $\omega \in \Omega$, an infinite-dimensional tuple whose $k$th element is $\omega_k \in \Omega_\iota$ which determines the stochastic gradient estimate.  The stochastic process defined by Algorithm~\ref{alg.sqp} can thus be expressed as
\bequationNN
  \{(X_k(\omega),G_k(\omega),H_k(\omega),D_k(\omega),D_k^\true(\omega),Y_k(\omega),Y_k^\true(\omega))\},
\eequationNN
where, for all $k \in \N{}$, the random variables are the iterate~$X_k(\omega)$, stochastic gradient estimator~$G_k(\omega)$, symmetric matrix $H_k(\omega)$, search direction~$D_k(\omega)$, true search direction $D_k^\true(\omega)$, Lagrange multiplier~$Y_k(\omega)$, and true Lagrange multiplier $Y_k^\true(\omega)$.  (As in Assumption~\ref{ass.tau_xi}, the true search direction and Lagrange multiplier are defined by the solution of \eqref{eq.system} if the stochastic gradient were replaced by the true gradient.) For instance, given $\omega$ that specifies a particular outcome of a run of the algorithm, the quantity $X_k(\omega) \in \R{n}$ is the $k$th iterate. Given initial conditions (including that $X_1(\omega) = x_1 \in \R{n}$ for all $\omega \in \Omega$), let ${\cal F}_1$ denote the $\sigma$-algebra corresponding to the initial conditions and, for all $k \in \N{} \setminus \{1\}$, let ${\cal F}_k$ denote the $\sigma$-algebra defined by the initial conditions and the random variables $\{G_1,\dots,G_{k-1}\}$.  We assume the following.

\begin{assumption}\label{ass.g}
  For all $k \in \N{}$, the stochastic gradient estimator satisfies $\E[G_k | {\cal F}_k] = \nabla f(X_k)$. In addition, there exists $\sigma \in \R{}_{>0}$ such that, for all $k \in \N{}$, it holds that $\E[\|G_k - \nabla f(X_k)\|_2^2 | \Fcal_k] \leq \sigma^2$.
\end{assumption}

\begin{assumption}\label{ass.H}
  For all $k \in \N{}$, the matrix $H_k \in \mathbb{S}^n$ is $\Fcal_k$-measurable and bounded in $\ell_2$-norm by $\kappa_H \in \R{}_{>0}$. In addition, there exists $\zeta \in (0, \kappa_H]$ such that, for all $k \in \N{}$, $u^TH_ku \geq \zeta \|u\|_2^2$ for all $u \in \Null(\nabla c(X_k)^T)$.
\end{assumption}

Under Assumptions~\ref{ass.g} and \ref{ass.H}, it follows that $\{\Fcal_k\}$ is a filtration for the probability space $(\Omega, \Fcal, \P)$.  In particular, the initial conditions and a realization of $\{G_1,\dots,G_{k-1}\}$ determine the realizations of $\{(X_j, H_j, D_j^\true, Y_j^\true)\}_{j=1}^k$ and $\{(D_j,Y_j)\}_{j=1}^{k-1}$, i.e., for all $k \in \N{}$, one has that $(X_k,H_k,D_k^\true,Y_k^\true)$ is $\Fcal_k$-measurable while $(G_k,D_k,Y_k)$ is $\Fcal_{k+1}$-measurable.  Examples of symmetric matrices satisfying Assumption~\ref{ass.H} for each $k \in \N{}$ are $H_k = I$ or $H_k$ being (an approximation of) the Hessian of the Lagrangian at $(X_k,Y_{k-1})$ as long as safeguards are included to ensure that $H_k$ is sufficiently positive definite in the null space of $\nabla c(X_k)^T$, as required by Assumption~\ref{ass.H}.

The following lemma characterizes decreases of the merit function~$\phi_\tau$ and the boundedness of $\Delta q_\tau$; it is adapted from Lemmas 3.7, 3.8, and 3.12 in~\cite{BCRZ21}.

\blemma\label{lem.sufficient_decrease}
  For all $k \in \N{}$, it follows that $\E[D_k | \Fcal_k] = D_k^\true$, 
  \begin{align}
    &\ \phi_\tau (X_k + \alpha_k D_k) - \phi_\tau(X_k) \nonumber \\
    \leq&\ -\alpha_k \Delta q_\tau(X_k, \nabla f(X_k), H_k, D_k^\true) \nonumber \\
    &\ + \thalf \alpha_k \beta_k \Delta q_\tau(X_k, G_k, H_k, D_k)
  + \alpha_k \tau \nabla f(X_k)^T(D_k - D_k^\true), \label{eq.merit_decr}
  \end{align}
  and
  \bequation \label{eq.merit_reduction_expec}
    \E[\Delta q_\tau (X_k, G_k, H_k, D_k) | \Fcal_k] \leq \Delta q_\tau (X_k, \nabla f(X_k), H_k, D_k^\true) + \tau \zeta^{-1} \sigma^2.
  \eequation
\elemma

\section{Convergence of primal iterates}\label{sec.primal_iterates}

One can derive from the analysis in \cite{BCRZ21} that, under Assumptions~\ref{ass.opt}, \ref{ass.tau_xi}, \ref{ass.g}, and \ref{ass.H}, Algorithm~\ref{alg.sqp} yields (recall the stationarity conditions \eqref{eq.kkt})
\bequation\label{eq.liminf_BCRZ}
  \liminf_{k\to\infty} \E[\|\nabla f(X_k) + \nabla c(X_k) Y_k^\true\|_2^2 + \|c(X_k)\|_2] = 0;
\eequation
see \cite[Corollary~3.14]{BCRZ21} for further details.  This lower limit in \eqref{eq.liminf_BCRZ} does not imply that the primal iterate sequence converges in any particular sense to a primal stationary point.  In this section, we establish conditions under which almost-sure convergence of the primal iterate sequence can be guaranteed, which when coupled with the results of the next section (in particular, Corollary~\ref{cor.y_conv}) provides a rich picture of the convergence behavior of Algorithm~\ref{alg.sqp}.

Before presenting the results of this section, let us state a common decomposition of the search direction that is used in our analysis.  Observe by the Fundamental Theorem of Linear Algebra that $D_k$ can be decomposed as $D_k = U_k + V_k$, where $U_k \in \Null(\nabla c(X_k)^T)$ and $V_k \in \Range(\nabla c(X_k))$.  Under Assumption~\ref{ass.opt}, we can now introduce $Z_k \in \R{n \times (n-m)}$ as a matrix whose columns form a basis for $\Null(\nabla c(X_k)^T)$.  It follows that $U_k = Z_k W_k$ for some $W_k \in \R{n-m}$.  Now, from~\eqref{eq.system}, it follows that
\bequationNN
  -c(X_k) = \nabla c(X_k)^TD_k = \nabla c(X_k)^T(Z_kW_k+V_k) = \nabla c(X_k)^T V_k.
\eequationNN
Since $\nabla c(X_k)^T$ has full row rank under Assumption~\ref{ass.opt}, the unique $V_k \in \Range(\nabla c(X_k))$ that solves this linear system is given by
\bequation\label{eq.v}
  V_k = -\nabla c(X_k) (\nabla c(X_k)^T \nabla c(X_k))^{-1} c(X_k) = -(\nabla c(X_k)^\dagger)^T c(X_k).
\eequation
Hence, from~\eqref{eq.system} and \eqref{eq.v}, it follows under Assumption~\ref{ass.H} that
\begin{align}
  && -G_k &= H_k(Z_kW_k+V_k) + \nabla c(X_k) Y_k \nonumber \\
  \implies && W_k &= - (Z_k^T H_k Z_k)^{-1} Z_k^T (G_k+H_kV_k) \nonumber \\
  \implies && Z_kW_k &= -Z_k(Z_k^TH_kZ_k)^{-1} Z_k^T (G_k+H_kV_k) \nonumber \\
  && &= -Z_k(Z_k^TH_kZ_k)^{-1} Z_k^T (G_k - H_k (\nabla c(X_k)^\dagger)^T c(X_k)). \label{eq.zw}
\end{align}

Our main result of this section is Theorem~\ref{thm.phi_conv} on page~\pageref{thm.phi_conv}.  Leading up to that result, we require two technical lemmas that are presented next.  The first, Lemma~\ref{lem.merit_conv} below, establishes a result about the asymptotic behavior of the sequence of merit function values, namely, $\{\phi_\tau(X_k)\}$, and about the sequence of reductions in the model of the merit function corresponding to the true objective gradients and corresponding search directions.  Note that the proof of the lemma uses the fact that the sequence $\{\beta_k\}$ employed in Algorithm~\ref{alg.sqp} is unsummable (i.e., $\sum_{k=1}^\infty \beta_k = \infty$), but square-summable (i.e., $\sum_{k=1}^\infty \beta_k^2 < \infty$).

\blemma \label{lem.merit_conv}
  Suppose that Assumptions~\ref{ass.opt}, \ref{ass.tau_xi}, \ref{ass.g}, and \ref{ass.H} hold.  Then,
  \begin{align*}
    \lim_{k \to \infty} \{\phi_\tau (X_k)\}\ \text{exists and is finite almost surely,} & \ \text{and} \\
    \liminf_{k \to \infty} \Delta q_\tau (X_k, \nabla f(X_k), H_k, D_k^\true ) = 0\  \text{almost surely.} &
  \end{align*}
\elemma
\proof
  For all $k \in \N{}$, define $R_k := \phi_\tau(X_k) - \tau f_{\inf}$, and observe from Assumption~\ref{ass.opt} that $R_k \geq 0$ for all $k \in \N{}$.  Also, for arbitrary $k \in \N{}$, one finds
  \begin{align}
     &\ \E[R_{k+1}| \Fcal_k] \nonumber \\
    =&\ \E[\phi_\tau(X_k + \alpha_k D_k) | \Fcal_k] - \tau f_{\inf} \nonumber \\
    \leq&\ \phi_\tau(X_k) - \tau f_{\inf} - \alpha_k \E[\Delta q_\tau(X_k, \nabla f(X_k), H_k, D_k^\true) | \Fcal_k] \nonumber \\
    &\ + \tfrac12 \alpha_k \beta_k \E[\Delta q_\tau(X_k, G_k, H_k, D_k) | \Fcal_k] \nonumber \\
    =&\ R_k - \alpha_k \Delta q_\tau (X_k, \nabla f(X_k), H_k, D_k^\true) + \tfrac12 \alpha_k \beta_k \E[\Delta q_\tau(X_k, G_k, H_k, D_k)|\Fcal_k] \nonumber \\
    \leq&\ R_k - \alpha_k \Delta q_\tau (X_k, \nabla f(X_k), H_k, D_k^\true) \nonumber \\
    &\ + \tfrac12 \alpha_k \beta_k \left( \Delta q_\tau (X_k, \nabla f(X_k), H_k, D_k^\true) + \tau \zeta^{-1} \sigma^2  \right) \nonumber \\
    =&\ R_k - \alpha_k (1 - \tfrac12 \beta_k) \Delta q_\tau (X_k, \nabla f(X_k), H_k, D_k^\true) + \tfrac12 \alpha_k \beta_k \tau \zeta^{-1} \sigma^2 \nonumber \\
    \leq&\ R_k - \alpha_k \Delta q_\tau (X_k, \nabla f(X_k), H_k, D_k^\true) + \tfrac12 \alpha_k \beta_k \tau \zeta^{-1} \sigma^2 \nonumber \\
    =&\ R_k - \tfrac{\beta_k \tau \xi}{\tau L_{\nabla f} + \Gamma} \Delta q_\tau (X_k, \nabla f(X_k), H_k, D_k^\true) + \tfrac{\beta_k^2 \tau^2 \xi \sigma^2}{2\zeta (\tau L_{\nabla f} + \Gamma)}, \label{eq.conv_z_bnd}
  \end{align}
  where the first inequality follows from Lemma~\ref{lem.sufficient_decrease}, the second follows from~\eqref{eq.merit_reduction_expec}, and the third follows from \eqref{eq.merit_reduction_lbnd} and the fact that $\{\beta_k\} \subset (0,1]$.  Now define
  \bequationNN
    P_k := \tfrac{\beta_k \tau \xi}{\tau L_{\nabla f} + \Gamma} \Delta q_\tau (X_k, \nabla f(X_k), H_k, D_k^\true)\ \ \text{and}\ \ Q_k := \tfrac{\beta_k^2 \tau^2 \xi \sigma^2}{2\zeta (\tau L_{\nabla f} + \Gamma)}
  \eequationNN
  and observe that $P_k \geq 0$ for all $k \in \N{}$ follows by~\eqref{eq.merit_reduction_lbnd}, $Q_k \geq 0$ for all $k \in \N{}$ follows since $\{\beta_k\} \subset (0,1]$, and $\sum_{k=1}^\infty Q_k < \infty$ since $\sum_{k=1}^\infty \beta_k^2 < \infty$.  Therefore, with \eqref{eq.conv_z_bnd}, applying Lemma~\ref{lem.rs} shows that, almost surely,
  \bequationNN
    \sum_{k=1}^\infty P_k = \sum_{k=1}^\infty \alpha_k \Delta q_\tau (X_k, \nabla f(X_k), H_k, D_k^\true) < \infty
  \eequationNN
  and $\lim_{k \to \infty} R_k = \lim_{k \to \infty} (\phi_\tau(X_k) - f_{\inf})$ exists and is finite.  The latter fact leads directly to the first desired conclusion, whereas the former fact along with $\sum_{k=1}^\infty \beta_k = \infty$ (so $\sum_{k=1}^\infty \alpha_k=\infty$) yields the second.
\qed

The components of our second technical lemma can be derived from various results from \cite{BCRZ21}, namely, Lemmas~2.10, 2.11, 2.12, and 3.4 in \cite{BCRZ21}.  Nonetheless, we provide a proof of the lemma for the sake of completeness.

\blemma\label{lem.d}
  Suppose that Assumptions~\ref{ass.opt}, \ref{ass.tau_xi}, \ref{ass.g}, and \ref{ass.H} hold.  Then, with respect to the decomposition $D_k^\true = U_k^\true + V_k$ for all $k \in \N{}$, where $U_k^\true \in \Null(\nabla c(X_k)^T)$ and $V_k \in \Range(\nabla c(X_k))$, and with
  \bequationNN
    \Psi_k := \begin{cases} \|U_k^\true\|_2^2 + \|c(X_k)\|_2 & \text{if $\|U_k^\true\|_2^2 \geq \kappa_{uv} \|V_k\|_2^2$} \\ \|c(X_k)\|_2 & \text{otherwise.} \end{cases}
  \eequationNN
  the following statements hold.
  \begin{enumerate}[label=(\alph*)]
    \item There exists $\kappa_{uv} \in \R{}_{>0}$ such that, for all $k \in \N{}$,
    \bequationNN
      \|U_k^\true\|_2^2 \geq \kappa_{uv} \|V_k\|_2^2 \implies (D_k^\true)^T H_k D_k^\true \geq \thalf \zeta \|U_k^\true\|_2^2.
    \eequationNN
    \item There exists $\kappa_\Psi \in \R{}_{>0}$ such that, for all $k \in \N{}$,
    \bequationNN
      \|D_k^\true\|_2^2 + \|c(X_k)\|_2 \leq (\kappa_\Psi + 1) \Psi_k.
    \eequationNN
    \item There exists $\kappa_q \in \R{}_{>0}$ such that, for all $k \in \N{}$,
    \bequationNN
      \Delta q_\tau(X_k, \nabla f(X_k), H_k, D_k^\true) \geq \kappa_q \Psi_k.
    \eequationNN
  \end{enumerate}
\elemma
\proof
  We prove each part in turn.
  \begin{enumerate}[label=(\alph*)]
    \item For arbitrary $k \in \N{}$ and $\kappa \in \R{}_{>0}$, one finds under Assumption~\ref{ass.H} that $\|U_k^\true\|_2^2 \geq \kappa \|V_k\|_2^2$ implies
    \begin{align*}
      (D_k^\true)^T H_k D_k^\true
      &= (U_k^\true)^T H_k U_k^\true + 2(U_k^\true)^T H_k V_k + V_k^T H_k V_k \\
      &\geq \zeta \|U_k^\true\|_2^2 - 2\|U_k^\true\|_2 \|H_k\|_2 \|V_k\|_2 - \|H_k\|_2 \|V_k\|_2^2 \\
      &\geq (\zeta - \tfrac{2\kappa_H}{\sqrt{\kappa}} - \tfrac{\kappa_H}{\kappa}) \|U_k^\true\|_2^2.
    \end{align*}
    Thus, the desired result holds for any $\kappa_{uv} \in \R{}_{>0}$ such that $\tfrac{2\kappa_H}{\sqrt{\kappa_{uv}}} + \tfrac{\kappa_H}{\kappa_{uv}} \leq \tfrac{\zeta}{2}$.
    \item If $\|U_k^\true\|_2^2 \geq \kappa_{uv} \|V_k\|_2^2$, it follows that
    \begin{align*}
      \|D_k^\true\|_2^2 = \|U_k^\true\|_2^2 + \|V_k\|_2^2
      &\leq (1+\kappa_{uv}^{-1}) \|U_k^\true\|_2^2 \\
      &\leq (1+\kappa_{uv}^{-1}) (\|U_k^\true\|_2^2 + \|c(X_k)\|_2).
    \end{align*}
    Otherwise, with \eqref{eq.v} and Assumption~\ref{ass.opt}, it follows that
    \begin{align*}
      \|D_k^\true\|_2^2 = \|U_k^\true\|_2^2 + \|V_k\|_2^2
      &< (\kappa_{uv}+1) \|V_k\|_2^2 \\
      &= (\kappa_{uv}+1) \|(\nabla c(X_k)^\dagger)^T c(X_k)\|_2^2 \\
      &\leq \kappa_c (\kappa_{uv}+1) r^{-2} \|c(X_k)\|_2.
    \end{align*}
    Hence, $\|D_k^\true\|_2^2 \leq \kappa_\Psi \Psi_k$ holds with $\kappa_\Psi := \max\{1 + \kappa_{uv}^{-1}, \kappa_c (\kappa_{uv}+1) r^{-2}\}$.  The desired conclusion then follows since, by definition, $\Psi_k \geq \|c(X_k)\|_2$.
    \item By~\eqref{eq.merit_reduction_lbnd} and part (a), it follows that $\|U_k^\true\|_2^2 \geq \kappa_{uv} \|V_k\|_2^2$ implies
    \bequationNN
      \Delta q_\tau (X_k, \nabla f(X_k), H_k, D_k^\true) \geq \tfrac14 \tau \zeta \|U_k^\true\|_2^2 + \nu \|c(X_k)\|_1,
    \eequationNN
    and otherwise one still finds $\Delta q_\tau (X_k, \nabla f(X_k), H_k, D_k^\true) \geq \nu \|c(X_k)\|_1$.  Hence, since $\|\cdot\|_2 \leq \|\cdot\|_1$, the result follows with $\kappa_q := \min\{\tfrac14 \tau \zeta, \nu \}$.
\qed
\end{enumerate}

We are now prepared to prove our main result of this section, which establishes conditions under which $\{\phi_\tau(X_k)\}$ converges almost-surely to a minimizer of the merit function at which the constraints of problem~\eqref{prob.opt} are satisfied and that $\{X_k\}$ converges almost-surely to a primal stationary point (in fact, a local minimizer).  After proving the result, we provide some additional commentary on the inequality~\eqref{eq.phi_ass} that is required for the theorem.

\btheorem\label{thm.phi_conv}
  Suppose that Assumptions~\ref{ass.opt}, \ref{ass.tau_xi}, \ref{ass.g}, and \ref{ass.H} hold.  In addition, suppose that there exists $x_\star \in \Xcal$ with $c(x_\star) = 0$, $\mu \in \R{}_{>1}$, and $\epsilon \in \R{}_{>0}$ such that
  for all $x \in \Xcal_{\epsilon,x_\star} := \{x \in \Xcal : \|x - x_\star\|_2 \leq \epsilon\}$ one finds
  \bequation\label{eq.phi_ass}
    \phi_\tau(x) - \phi_\tau(x_\star) \bcases = 0 & \text{if $x = x_\star$} \\ \in (0, \mu (\tau \|Z(x)^T \nabla f(x)\|_2^2 + \|c(x)\|_2)] & \text{otherwise,} \ecases
  \eequation
  where for all $x \in \Xcal_{\epsilon,x_\star}$ one defines $Z(x) \in \R{n \times (n-m)}$ as some orthonormal matrix whose columns form a basis for the null space of $\nabla c(x)^T$.  Then, if $\displaystyle \limsup_{k\to\infty} \{\|X_k - x_\star\|_2\} \leq \epsilon$ almost surely, it follows that
  \begin{align*}
    \{\phi_\tau (X_k)\} \asto \phi_\tau(x_\star),\ \ \{X_k\} &\asto x_\star, \\
    \text{and}\ \ \left\{ \bbmatrix \nabla f(X_k) + \nabla c(X_k)Y_k^\true \\ c(X_k) \ebmatrix \right\} &\asto 0.
  \end{align*}
\etheorem
\proof
  As before, for all $k \in \N{}$, let $Z_k \in \R{n \times (n-m)}$ denote an orthonormal matrix whose columns form a basis for $\Null(\nabla c(X_k))$, where if $X_k \in \Xcal_{\epsilon,x_\star}$ then $Z_k = Z(X_k)$ is one such that \eqref{eq.phi_ass} holds at $x = X_k$.  For arbitrary $k \in \N{}$, it follows from~\eqref{eq.system} that the reduced gradient satisfies
  \bequationNN
    Z_k^T \nabla f(X_k) = -Z_k^TH_kU_k^\true - Z_k^T H_k V_k,
  \eequationNN
  from which it follows with \eqref{eq.v}, \eqref{eq.zw}, and Assumptions~\ref{ass.opt} and \ref{ass.H} that
  \begin{align*}
    \|Z_k^T \nabla f(X_k)\|_2^2
    &= \|Z_k^TH_kU_k^\true + Z_k^T H_k V_k\|_2^2 \\
    &\leq 2(\|Z_k^TH_kU_k^\true\|_2^2 + \|Z_k^T H_k V_k\|_2^2) \\
    &\leq 2 \kappa_H^2 (\|U_k^\true\|_2^2 + \|V_k\|_2^2) \\
    &\leq 2 \kappa_H^2
      (\|U_k^\true\|_2^2 + \|(\nabla c(X_k)^\dagger)^T\|_2^2 \|c(X_k)\|_2^2) \\
    &\leq 2 \kappa_H^2
      (\|U_k^\true\|_2^2 + \kappa_c r^{-2} \|c(X_k)\|_2).
  \end{align*}
  Combining this with \eqref{eq.phi_ass} and Lemma~\ref{lem.d}, one finds that $X_k \in \Xcal_{\epsilon,x_\star}$ implies
  \begin{align*}
    &\ \phi_\tau(X_k) - \phi_\tau(x_\star) \\
    \leq&\ \mu (2 \tau \kappa_H^2 \|U_k^\true\|_2^2 + (2 \tau \kappa_c \kappa_H^2 r^{-2} + 1) \|c(X_k)\|_2) \\
    \leq&\ \mu \max\{2 \tau \kappa_H^2, 2 \tau \kappa_c \kappa_H^2 r^{-2} + 1\} ( \|D_k^\true\|_2^2 + \|c(X_k)\|_2) \\
    \leq&\ \mu \max\{2 \tau \kappa_H^2, 2 \tau \kappa_c \kappa_H^2 r^{-2} + 1\} (\kappa_\Psi + 1) \Psi_k \\
    \leq&\ \mu \max\{2 \tau \kappa_H^2, 2 \tau \kappa_c \kappa_H^2 r^{-2} + 1\} (\kappa_\Psi + 1) \kappa_q^{-1} \Delta q_\tau(X_k, \nabla f(X_k), H_k, D_k^\true).
  \end{align*}
  Since by Lemma~\ref{lem.merit_conv} one has $\displaystyle \liminf_{k\to\infty} \Delta q_\tau(X_k, \nabla f(X_k), H_k, D_k^\true) = 0$ almost surely, it follows from above and the conditions of the theorem that the limit $\displaystyle \lim_{k \to \infty} \phi_\tau (X_k)$ (which exists almost surely by Lemma~\ref{lem.merit_conv}) must be $\phi_\tau(x_\star)$, as desired.  The remaining conclusions follow from $\{\phi_\tau (X_k)\} \asto \phi_\tau(x_\star)$ and the facts that, by \eqref{eq.phi_ass}, the point $x_\star$ is the unique point in $\Xcal_{\epsilon,x_\star}$ with merit function value equal to $\phi_\tau(x_\star)$ and that is a primal stationary point.
\qed

Observe that \eqref{eq.phi_ass} can be viewed as a generalization of the well-known Polyak--Łojasiewicz condition from the unconstrained continuous optimization literature \cite{Poly63}.  Indeed, if the constraints are affine and one considers $x$ such that $x - x_\star \in \Null(\nabla c(x_\star)^T)$, then~\eqref{eq.phi_ass} says that the squared $\ell_2$-norm of the reduced gradient $Z(x)^T\nabla f(x)$ is at least proportional to $\phi_\tau(x) - \phi_\tau(x_\star)$.  On the other hand, if the objective function remains constant along displacements in $\Range(\nabla c(x))$, then \eqref{eq.phi_ass} says that the norm of the constraint violation is at least proportional to this difference in merit function values.  More generally, such as for nonlinear $f$ and $c$, the condition \eqref{eq.phi_ass} is a generalization of these special cases that says that the combination of the squared $\ell_2$-norm of the reduced gradient and constraint violation is at least proportional to the optimality gap in the merit function in a neighborhood of the point $x_\star$.

\section{Convergence of Lagrange multipliers}\label{sec.multipliers}

In this section, we study convergence properties of the sequence $\{Y_k\}$ generated by Algorithm~\ref{alg.sqp}.  Let us begin by expressing the solution component $Y_k$ of \eqref{eq.system} for arbitrary $k \in \N{}$ using the step decomposition stated in \eqref{eq.v} and \eqref{eq.zw}.  One finds by substituting~\eqref{eq.v} and \eqref{eq.zw} back into~\eqref{eq.system} gives
\begin{align*}
  \nabla c(X_k) Y_k
  &= -(H_k(Z_kW_k+V_k)+G_k) \\
  &= (I-H_k Z_k (Z_k^T H_k Z_k)^{-1} Z_k^T) (H_k (\nabla c(X_k)^\dagger)^T c(X_k) - G_k),
\end{align*}
from which it follows under Assumption~\ref{ass.opt} that
\bequation\label{eq.y_prelim}
  Y_k = \nabla c(X_k)^\dagger (I-H_k Z_k (Z_k^T H_k Z_k)^{-1} Z_k^T) (H_k (\nabla c(X_k)^\dagger)^T c(X_k) - G_k).
\eequation
For notational convenience, let us now define the matrix
\bequation\label{eq.M}
  M_k := \nabla c(X_k)^\dagger (I - H_k Z_k (Z_k^T H_k Z_k)^{-1} Z_k^T) \in \R{m \times n},
\eequation
which can be viewed as the product of a pseudoinverse and a projection matrix, so that we may succinctly write from \eqref{eq.y_prelim} that
\bequation\label{eq.y}
  Y_k = M_k (H_k (\nabla c(X_k)^\dagger)^T c(X_k) - G_k).
\eequation
(We remark in passing that an alternative to our subsequent discussions and analysis could be considered where, for all $k \in \N{}$, the vector $Y_k$ is not defined through \eqref{eq.system}, but rather is set as a so-called least squares multiplier, i.e., as $Y_k = \arg\min_{Y\in\R{m}} \|G_k + \nabla c(X_k)Y\|_2^2$.  For one thing, this choice removes the dependence of $Y_k$ on $H_k$, which might have certain advantages.  However, for our purposes, we focus on multipliers being computed through \eqref{eq.system} since this is a popular approach in practice and does not require additional computation.)

Our goal in our analysis of properties of Lagrange multiplier estimators is to prove convergence of such estimators when the primal iterate approaches a primal stationary point for~\eqref{prob.opt}.  After all, it is only when the primal iterate lies in such a neighborhood (or at least a neighborhood of the feasible region) that the Lagrange multiplier has meaning in terms of certifying stationarity.  Consequently, a main focus in our results is how~$Y_k$ may be viewed through~\eqref{eq.y} as a function of the primal iterate $X_k$.  For our analysis here, we make the following assumption, where given $x \in \R{n}$ and $\epsilon \in \R{}_{>0}$ we define (consistent with Theorem~\ref{thm.phi_conv}) the neighborhood $\Xcal_{\epsilon,x} := \{\xbar \in \R{n} : \|\xbar - x\|_2 \leq \epsilon\}$.

\bassumption\label{ass.M}
  Given $x_\star \in \Xcal$ as a primal stationary point for problem~\eqref{prob.opt}, there exist $\epsilon \in \R{}_{>0}$, $\Hcal : \R{n} \to \mathbb{S}^n$, $L_\Hcal \in \R{}_{>0}$, $\Mcal : \R{n} \to \R{m \times n}$, and $L_\Mcal \in \R{}_{>0}$ such that the following hold.
  \benumerate
    \item[(i)] $H_k = \Hcal(X_k)$ whenever $X_k \in \Xcal_{\epsilon,x_\star}$;
    \item[(ii)] $\|\Hcal(x) - \Hcal(\xbar)\|_2 \leq L_\Hcal \|x - \xbar\|_2$ for all $(x,\xbar) \in \Xcal_{\epsilon,x_\star} \times \Xcal_{\epsilon,x_\star}$;
    \item[(iii)] $M_k = \Mcal(X_k)$ whenever $X_k \in \Xcal_{\epsilon,x_\star}$; and
    \item[(iv)] $\|\Mcal(x) - \Mcal(\xbar)\|_2 \leq L_\Mcal \|x - \xbar\|_2$ for all $(x,\xbar) \in \Xcal_{\epsilon,x_\star} \times \Xcal_{\epsilon,x_\star}$.
  \eenumerate
\eassumption

A few important observations and justifications are in order.

\bitemize
  \item Assumption~\ref{ass.M} states that, in a neighborhood of a primal stationary point~$x_\star$, the algorithm sets the matrix $H_k$ through a Lipschitz continuous function of $X_k$.  This clearly holds if the algorithm chooses $H_k = H$ for all $k \in \N{}$ for some prescribed $H \in \mathbb{S}_{\succ0}^n$ that guarantees that Assumption~\ref{ass.H} holds for all $k \in \N{}$.  Alternatively, the requirements on $\{H_k\}$ in Assumption~\ref{ass.M} can be satisfied if second-order derivatives of $f$ are Lipschitz continuous and the algorithm chooses $H_k$ as (an approximation of) the Hessian of the objective at $X_k$, or even as (an approximation of) the Hessian of the Lagrangian at $(X_k,\Ybar_k)$ if the second-order derivatives of the components of $c$ are also Lipschitz continuous and care is taken in the selection of $\Ybar_k$.  In practice, one might consider $\Ybar_k = Y_{k-1}$ for all $k \in \N{}$, but as one can see, this choice is influenced by noise in the stochastic gradient estimators, which can cause $\{H_k\}$ to violate Assumption~\ref{ass.H}.  An alternative choice that would not violate the assumption is to choose $\Ybar_k$ as a prescribed vector or at least one that remains fixed for all sufficiently large $k \in \N{}$.  In any case, one finds that there exist reasonable choices for the algorithm to choose $\{H_k\}$ such that Assumption~\ref{ass.M} holds.  (We remark in passing that deterministic Newton-based methods for solving constrained optimization problems can possess local convergence guarantees of the Lagrange multipliers under relatively strong assumptions when $H_k$ is the Hessian of the Lagrangian at $(X_k,Y_{k-1})$, and this might also be achievable for a stochastic algorithm with highly accurate gradient estimates.  However, since our algorithm of interest operates in a highly stochastic regime in which only Assumption~\ref{ass.g} is assumed to hold, allowing unadulterated use of the Hessian of the Lagrangian at $(X_k,Y_{k-1})$ is not reasonable.)
  \item Assumption~\ref{ass.M} also states that, in a neighborhood of a primal stationary point~$x_\star$, the matrix $M_k$ is defined through a Lipschitz continuous function of $X_k$. Given the assumptions about $\Hcal$, Assumption~\ref{ass.opt} (which states that the constraint Jacobians have full row rank and the constraint Jacobian function is Lipschitz continuous over the set~$\Xcal$ containing the iterates), and Lemma~\ref{lem.pinv}, it might seem at first glance that this assumption about~$\Mcal$ is a straightforward consequence of these other assumptions.  However, one needs to be careful.  To understand a potential issue when trying to draw such a conclusion, recall that a product of functions that are Lipschitz continuous on a bounded set is itself Lipschitz continuous on the set.  Hence, one might attempt to justify our Lipschitz-continuity assumption about $\Mcal$ by assuming that, for all $k \in \N{}$, the null space basis $Z_k$ can be viewed as being generated by a Lipschitz continuous function of $X_k$.  Note, however, that it has been shown in \cite{BS86} that, even when $\nabla c(\cdot)^T$ is continuous with values that have full row rank for all $x \in \R{n}$, there does not necessarily exist continuous $\Zcal : \R{n} \to \R{n \times (n-m)}$ such that for all $x \in \R{n}$ the columns of $\Zcal(x)$ form a basis for $\Null(\nabla c(x)^T)$.  That being said, as also shown in \cite{BS86}, it is possible to employ procedures such that the reduced Hessian (approximation), namely, $Z_k^TH_kZ_k$, does not depend on the choice of $Z_k$, and, given a point $\xbar \in \R{n}$ at which $\nabla c(\xbar)^T$ has full row rank, there exists a neighborhood of $\xbar$ over which one can define a ``null space basis function'' that is continuous over the neighborhood.  See also \cite{BN91}, which discusses two procedures that ensure that, in the neighborhood of a given point, a null space basis can be defined by a Lipschitz continuous null space basis function.  Overall, due to these observations, we contend that assuming Lipschitz continuity of~$\Mcal$---in a sufficiently small neighborhood of $x_\star$, as is stated in Assumption~\ref{ass.M}---is justified for our analysis.
\eitemize

We are now prepared to prove that, near a primary stationary point, the expected error in the Lagrange multiplier estimator is bounded by the distance of the primal iterate to the primal stationary point plus an error due to the stochastic gradient estimator.  The following theorem may be viewed as the main result of this section since the subsequent results follow from it.

\btheorem\label{thm.multiplier}
  Suppose that Assumptions~\ref{ass.opt}, \ref{ass.tau_xi}, \ref{ass.g}, \ref{ass.H}, and \ref{ass.M} hold, $(x_\star,y_\star)$ is a stationary point for problem~\eqref{prob.opt}, and $\epsilon \in \R{}_{>0}$ is defined as in Assumption~\ref{ass.M}.  Then, for any $k \in \N{}$, one finds $\|X_k - x_\star\|_2 \leq \epsilon$ implies
  \bequation \label{eq.y_conv_ineq}
    \|Y_k-y_\star\|_2 \leq \kappa_y \|X_k - x_\star\|_2 + r^{-1} \|\nabla f(X_k) - G_k\|_2,
  \eequation
  where $\kappa_y := \kappa_H L_c r^{-2} + L_{\nabla f} r^{-1} + \kappa_{\nabla f} L_\Mcal$.
\etheorem
\proof
  Under the conditions of the theorem, $\|X_k - x_\star\|_2 \leq \epsilon$ and \eqref{eq.y} imply that $Y_k = M_k(H_k(\nabla c(X_k)^\dagger)^Tc(X_k) - G_k)$, where $M_k = \Mcal(X_k)$ and $H_k = \Hcal(X_k)$, while one can similarly derive that $y_\star = -M_\star \nabla f(x_\star)$, where $M_\star = \Mcal(x_\star)$. Hence, the difference $Y_k - y_\star$ can be decomposed into three parts as follows:
  \bequation\label{eq.y_decomp}
    \baligned
      Y_k-y_\star =&\ M_k H_k (\nabla c(X_k)^\dagger)^T c(X_k) \\
      &\ + M_k (\nabla f(x_\star) - G_k) + (M_\star - M_k) \nabla f(x_\star).
    \ealigned
  \eequation
  Consequently, to bound $\|Y_k-y_\star\|_2$ when $\|X_k - x_\star\|_2 \leq \epsilon$, one can employ the triangle inequality and bound the norms of the three terms on the right-hand side of~\eqref{eq.y_decomp} separately. First, one finds that
  \begin{align}
    &\ \|M_k H_k (\nabla c(X_k)^\dagger)^T c(X_k)\|_2 \nonumber \\
    \leq&\ \|\nabla c(X_k)^\dagger\|_2 \|I - H_k Z_k (Z_k^TH_kZ_k)^{-1} Z_k^T\|_2 \|H_k\|_2 \|(\nabla c(X_k)^\dagger)^T\|_2 \|c(X_k)\|_2 \nonumber \\
    \leq&\ \kappa_H r^{-2} \|c(X_k) - c(x_\star)\|_2 \leq \kappa_H L_c r^{-2} \|X_k-x_\star\|_2, \label{eq.y_conv_prf_1}
  \end{align}
  where the first inequality follows from \eqref{eq.M} and properties of norms; the second follows from $c(x_\star) = 0$, Assumptions~\ref{ass.opt} and \ref{ass.H}, and the fact that $P_k := I-H_kZ_k(Z_k^TH_kZ_k)^{-1}Z_k^T$ is a projection matrix, so $\|P_k\|_2 \leq 1$; and the last follows from Lipschitz continuity of $c$ (see~\eqref{eq.Lipschitz}).  Second,
  \begin{align}
    &\ \|M_k (\nabla f(x_\star) - G_k)\|_2 \nonumber \\
    \leq&\ \|\nabla c(X_k)^\dagger\|_2 \|I-H_k Z_k (Z_k^TH_kZ_k)^{-1} Z_k^T\|_2 \|\nabla f(x_\star) - G_k \|_2 \nonumber \\
    \leq&\ r^{-1} \|\nabla f(x_\star) - G_k\|_2 \nonumber \\
    \leq&\ r^{-1} (\|\nabla f(x_\star)-\nabla f(X_k)\|_2 + \|\nabla f(X_k) - G_k\|_2) \nonumber \\
    \leq&\ L_{\nabla f} r^{-1} \|X_k-x_\star\|_2 + r^{-1} \|\nabla f(X_k) - G_k\|_2, \label{eq.y_conv_prf_2}
  \end{align}
  where the first inequality follows from \eqref{eq.M} and properties of norms, the second follows from Assumption~\ref{ass.opt} and $\|P_k\|_2 \leq 1$ (see above), the third follows from $\nabla f(x_\star)-G_k = \nabla f(x_\star) - \nabla f(X_k) + \nabla f(X_k) - G_k$ and the triangle inequality, and the last follows from Lipschitz continuity of $\nabla f$ (see~\eqref{eq.Lipschitz}).  Third,
  \bequation\label{eq.y_conv_prf_3}
    \|(M_\star - M_k) \nabla f(x_\star)\|_2 \leq \|M_\star - M_k\|_2 \|\nabla f(x_\star)\|_2 \leq \kappa_{\nabla f} L_\Mcal \|X_k - x_\star\|_2
  \eequation
  follows from properties of norms and Assumptions~\ref{ass.opt} and \ref{ass.M}. Substituting~\eqref{eq.y_conv_prf_1}, \eqref{eq.y_conv_prf_2}, and \eqref{eq.y_conv_prf_3} into~\eqref{eq.y_decomp} yields the desired conclusion.
\qed

Theorem~\ref{thm.multiplier} shows that whenever $X_k$ is sufficiently close to a primal stationary point $x_\star$, the Lagrange multiplier error $\|Y_k-y_\star\|_2$ is bounded by a constant times the primal iterate error $\|X_k-x_\star\|_2$ plus a constant times the error in the stochastic gradient estimator. The latter term is inevitable since, no matter the proximity of $X_k$ to $x_\star$, the Lagrange multiplier estimator $Y_k$ is influenced by~$G_k$.  To emphasize this point, we present the following corollary showing that, under the same conditions as Theorem~\ref{thm.multiplier}, the Lagrange multiplier estimator $Y_k^\true$---which appears in the convergence guarantee \cite[Corollary~3.14]{BCRZ21}, even though it is a quantity that is not actually computed in the algorithm---satisfies the same bound as $Y_k$, except without the additional noise.

\begin{corollary}\label{cor.y_conv_true}
  Suppose that Assumptions~\ref{ass.opt}, \ref{ass.tau_xi}, \ref{ass.g}, \ref{ass.H}, and \ref{ass.M} hold, $(x_\star,y_\star)$ is a stationary point for problem~\eqref{prob.opt}, and $\epsilon \in \R{}_{>0}$ is defined as in Assumption~\ref{ass.M}.  Then, for any $k \in \N{}$, one finds $\|X_k - x_\star\|_2 \leq \epsilon$ implies
  \bequation
  \|Y_k^\true - y_\star\|_2 \leq \kappa_y \|X_k - x_\star\|_2, \label{eq.y_conv_true}
  \eequation
  where $\kappa_y \in \R{}_{>0}$ is defined as in Theorem~\ref{thm.multiplier}.
\end{corollary}
\proof
  Following the same steps as in the proof of Theorem~\ref{thm.multiplier}, the difference $Y_k^\true - y_\star$ can now be decomposed as in~\eqref{eq.y_decomp} with $G_k = \nabla f(X_k)$. Therefore, in place of \eqref{eq.y_conv_prf_2}, one instead obtains the bound
  \begin{align*}
    &\ \|\nabla c(X_k)^\dagger (I-H_k Z_k (Z_k^T H_k Z_k)^{-1} Z_k^T) (\nabla f(x_\star) - \nabla f(X_k))\|_2 \\
    \leq&\ r^{-1} \|\nabla f(x_\star) - \nabla f(X_k)\|_2 \leq L_{\nabla f} r^{-1} \|X_k-x_\star\|_2,
  \end{align*}
  which when combined with~\eqref{eq.y_conv_prf_1}, \eqref{eq.y_conv_prf_3}, and~\eqref{eq.y_decomp} gives the desired result.
\qed

The previous theorem and corollary lead to the following corollary, which considers the special case that there exists an iteration number beyond which the primal iterates remain within a neighborhood of a particular primal stationary point almost surely.  In such a case, the expected error in the Lagrange multiplier is bounded by a constant times the distance from the primal iterate to the primal stationary point plus a term that is proportional to the bound on the noise of the stochastic gradient estimators.  The corollary also adds that the ``true'' multiplier does not have this additional error term.  (We remind the reader that, in Section~\ref{sec.primal_iterates}, we have shown conditions under which one can conclude that $\{X_k\} \asto x_\star$.  While this convergence guarantee does not directly imply that $\kbar \in \N{}$ exists as in the following corollary, it motivates consideration of this special case as one of practical relevance.)

\begin{corollary}\label{cor.y_conv_neighborhood}
  Suppose that Assumptions~\ref{ass.opt}, \ref{ass.tau_xi}, \ref{ass.g}, \ref{ass.H}, and \ref{ass.M} hold, $(x_\star,y_\star)$ is a stationary point for problem~\eqref{prob.opt}, $\epsilon \in \R{}_{>0}$ is defined as in Assumption~\ref{ass.M}, and for some $\kbar \in \N{}$ one has almost surely that
  \bequationNN
    \|X_k - x_\star\|_2 \leq \epsilon\ \text{for all}\ k \in \N{}\ \text{with}\ k \geq \kbar.
  \eequationNN
  Then, for any $k \in \N{}$ with $k \geq \kbar$, one finds that
  \begin{align*}
    \E[\|Y_k - y_\star\|_2 | \Fcal_k] &\leq \kappa_y \|X_k - x_\star\|_2 + r^{-1} \sigma \\ \text{and}\ \ 
    \E[\|Y_k^\true - y_\star\|_2 | \Fcal_k] &\leq \kappa_y \|X_k - x_\star\|_2, 
  \end{align*}
  where $\kappa_y \in \R{}_{>0}$ is defined as in Theorem~\ref{thm.multiplier}.
\end{corollary}
\proof
  Assumption~\ref{ass.g} and Jensen's inequality imply for all $k \in \N{}$ that
  \bequationNN
    \E[\|G_k -\nabla f(X_k)\|_2 | \Fcal_k]^2 \leq \E[\|G_k -\nabla f(X_k)\|_2^2 | \Fcal_k] \leq \sigma^2.
  \eequationNN
  The desired conclusions follow from this fact, the conditions of the corollary, and the conclusions of Theorem~\ref{thm.multiplier} and Corollary~\ref{cor.y_conv_true}.
\qed

At this point, we have determined useful properties of two multiplier estimator sequences, $\{Y_k\}$ and $\{Y_k^\true\}$.  In the neighborhood of a stationary point $(x_\star,y_\star)$, each element $Y_k^\true$ of the latter sequence has error that is proportional to the distance from $X_k$ to the primal stationary point $x_\star$, which is as much as can be expected.  However, this estimator is not realized by the algorithm as it requires explicit knowledge of the true gradient $\nabla f(X_k)$, which is not realized by the algorithm.  Each element $Y_k$ of the former sequence has the same error plus an additional error due to the stochastic gradient estimator.  One can reduce this error by computing a more accurate gradient estimate, say, using a large mini-batch of gradient estimates.  However, to avoid the need of computing a highly accurate gradient estimate in any given iteration, one might ask whether it is possible to obtain an estimator with lower error merely under Assumption~\ref{ass.g}?  Using classical ideas from statistical estimation, specifically mean estimation, a natural idea to consider is whether a more accurate estimator can be defined by an average, say, $Y_k^\avg := \tfrac{1}{k} \sum_{i=1}^k Y_k$.  After all, if the stochastic process $\{X_k\}$ converges in some sense to $x_\star$, then with $Y_k$ defined through a continuous function of $X_k$ in a neighborhood of $x_\star$ for all $k \in \N{}$, one might expect $\{Y_k^\avg\}$ to converge in the same sense to $y_\star$.

The main challenge in the analysis of $\{Y_k^\avg\}$ is that, in contrast to the classical setting of mean estimation using independent and identically distributed random variables, the multipliers $\{Y_k\}$ are computed at points $\{X_k\}$ at which the distributions of $\{G_k\}$ (or $\{\nabla f(X_k)\}$ for that matter) are neither independent nor identically distributed.  Hence, our analysis of the averaged sequence requires some additional assumptions that are introduced in the following lemma.  We prove the lemma and a following theorem, then discuss and provide justifications for these required assumptions in further detail.

\blemma\label{lem.noise_avg}
  Suppose that Assumptions~\ref{ass.opt}, \ref{ass.tau_xi}, \ref{ass.g}, \ref{ass.H}, and \ref{ass.M} hold, and with $M_k$ defined in \eqref{eq.M} and $\Delta_k := \nabla f(X_k) - G_k$ for all $k \in \N{}$ one finds
  \begin{subequations}\label{eq.noise_avg_ass}
    \begin{align}
      & \tfrac{1}{k} \E [ \|M_i \Delta_i\|_2^2 ] < \infty \text{ for all } (k, i) \in \N{} \times [k], \label{eq.yavg_ass_1} \\
      & \left\{\tfrac{1}{k} \sum_{i=1}^k \E\Bigg[ \|M_i \Delta_i\|_2^2 \mathbf{1}_{\left\{ \tfrac{\|M_i \Delta_i\|_2}{\sqrt k} > \delta \right\}} \Bigg| \Fcal_i \Bigg] \right\}\pto 0\ \text{for all}\ \delta \in \R{}_{>0}, \label{eq.yavg_ass_2} \\
      & \left\{\tfrac{1}{k} \sum_{i=1}^k \E [ M_i \Delta_i \Delta_i^T M_i^T | \Fcal_i ]\right\} \pto \Sigma\ \text{for some}\ \Sigma \in \mathbb{S}^n,\ \text{and} \label{eq.yavg_ass_3} \\
      & \sup_{k \in \N{}} \E \left[ \left\|\sum_{i=1}^k \tfrac{1}{\sqrt k} M_i \Delta_i\right\|_2^\Theta \right] < \infty\ \text{for some}\ \Theta \in \Rmbb_{>1}. \label{eq.yavg_ass_4}
    \end{align}
  \end{subequations}
  Then, it follows that
  \bequation\label{eq.noise_avg}
    \lim_{k \to \infty} \tfrac{1}{k} \E \left[ \left \|\sum_{i=1}^k M_i \Delta_i \right\|_2 \right] = 0.
  \eequation
\elemma
\proof
  For $(k,i) \in \N{} \times [k]$, let $\Fcal_{k,i} := \Fcal_i$ and $\xi_{k,i} := \tfrac{1}{\sqrt{k}} M_i \Delta_i$.  Under Assumption \ref{ass.g}, it follows that $\E[\xi_{k,i} | \Fcal_{k,i}] = 0$ for all $(k,i) \in \N{} \times [k]$ and that $\{(\xi_{k,i}, \Fcal_{k,i})\}_{k \in \N{},i \in [k]}$ is an $n$-dimensional martingale difference triangular array (see Lemma~\ref{lem.mart_clt}). Employing Lemma~\ref{lem.mart_clt}, one finds that 
  \bequationNN
    \left\{\sum_{i=1}^k \xi_{k,i}\right\} = \left\{\sum_{i=1}^k \tfrac{1}{\sqrt k} M_i \Delta_i\right\} \dto \Ncal(0, \Sigma).
  \eequationNN
  This fact, Lemma~\ref{lem.moment_convergence}, and \eqref{eq.yavg_ass_4} together yield
  \bequationNN
    \lim_{k \to \infty} \E \left[ \left\| \sum_{i=1}^k  \tfrac{1}{\sqrt k} M_i \Delta_i \right\|_2 \right] = \tilde \sigma < \infty, 
  \eequationNN
  where $\tilde \sigma := \Embb[\|X\|_2]$ for $X \sim \Ncal(0, \Sigma)$.  Therefore,
  \bequationNN
  \lim_{k \to \infty} \tfrac{1}{k} \E \left[ \left\| \sum_{i=1}^k  M_i \Delta_i \right\|_2 \right] = \(\lim_{k \to \infty} \tfrac{1}{\sqrt k} \) \( \lim_{k \to \infty}\E \left[ \left\| \tfrac{1}{\sqrt k}\sum_{i=1}^k  M_i \Delta_i \right\|_2 \right] \) = 0,
  \eequationNN
  which gives the desired result~\eqref{eq.noise_avg}.
\qed

The previous lemma allows us to prove the following theorem.  We remark upfront that the iteration index $\kbar \in \N{}$ that is used in the definition of $\{Y_k^\avg\}_{k=\kbar}^\infty$ in the theorem is not necessarily known in advance by the algorithm.  In particular, it is not necessarily known by the algorithm when the primal iterate sequence may be guaranteed to have entered a neighborhood of a primary stationary point such that the conditions of Assumption~\ref{ass.M} hold with a certain $\epsilon \in \R{}_{>0}$.  A stronger assumption would be to assume that $\kbar = 1$, i.e., that the algorithm is initialized at a point within a sufficiently small neighborhood of $x_\star$ in which it remains almost surely.  However, we state the result for general $\kbar$ for the sake of generality, and since in practice it would be reasonable to consider a scheme that only computes an averaged Lagrange multiplier estimate over those recent multipliers such that the corresponding primal iterates are within a prescribed neighborhood of the current iterate.  We consider such a strategy in our experiments in Section~\ref{sec.numerical}.

\btheorem\label{cor.y_conv_avg}
  Suppose that Assumptions~\ref{ass.opt}, \ref{ass.tau_xi}, \ref{ass.g}, \ref{ass.H}, and \ref{ass.M} hold, $(x_\star,y_\star)$ is a stationary point for problem~\eqref{prob.opt}, $\epsilon \in \R{}_{>0}$ is defined as in Assumption~\ref{ass.M}, for some $\kbar \in \N{}$ one has almost surely that
  \bequationNN
    \|X_k - x_\star\|_2 \leq \epsilon\ \text{for all}\ k \in \N{}\ \text{with}\ k \geq \kbar,
  \eequationNN
  $\{M_k\}$ and $\{\Delta_k\}$ are defined as in Lemma~\ref{lem.noise_avg}, and \eqref{eq.noise_avg_ass} holds.  Then, with $Y_k^\avg := \tfrac{1}{k-\kbar+1} \sum_{i=\kbar}^k Y_i$ for all $k \in \N{}$ with $k \geq \kbar$, one finds that
  \bequationNN
    \limsup_{k \geq \kbar, k \to \infty} \E[\|Y_k^\avg - y_\star\|_2] \leq \limsup_{k \geq \kbar, k \to \infty} \tfrac{\kappa_y}{k-\kbar+1} \sum_{i=\kbar}^k \E[\|X_i - x_\star\|_2].
  \eequationNN
\etheorem
\proof
  As in the proof of Theorem~\ref{thm.multiplier}, one finds for any $k \in \N{}$ with $k \geq \kbar$ that
  \begin{align}
    \|Y_k^\avg - y_\star\|_2
      =&\ \tfrac{1}{k-\kbar+1} \left\| \sum_{i=\kbar}^k (Y_i - y_\star) \right\|_2 \nonumber\\
      \leq&\ \tfrac{\kappa_H L_c r^{-2} + \kappa_{\nabla f} L_\Mcal}{k-\kbar+1} \sum_{i=\kbar}^k \|X_i - x_\star\|_2 \nonumber\\
      &\ + \tfrac{1}{k-\kbar+1} \left\|\sum_{i=\kbar}^k M_i (\nabla f(x_\star) - \nabla f(X_k) + \nabla f(X_k) - G_k) \right\|\nonumber \\
      \leq&\ \tfrac{\kappa_y}{k-\kbar+1} \sum_{i=\kbar}^k \|X_i - x_\star\|_2 + \tfrac{1}{k-\kbar+1} \left\|\sum_{i=\kbar}^k M_i \Delta_i \right\|_2 \label{eq.y_conv_avg_pf},
  \end{align}
  Thus, under the conditions of the theorem, it follows with Lemma~\ref{lem.noise_avg} that
  \begin{align*}
    &\ \limsup_{k\geq\kbar, k \to \infty} \E[\|Y_k^\avg - y_\star\|_2] \\
    \leq&\ \limsup_{k\geq\kbar, k \to \infty} \tfrac{\kappa_y}{k-\kbar+1} \sum_{i=\kbar}^k \E[\|X_i-x_\star\|_2] + \limsup_{k\geq\kbar, k \to \infty} \tfrac{1}{k-\kbar+1} \E \left[ \left\| \sum_{i=\kbar}^k M_i \Delta_i \right\|_2 \right] \\
    =&\ \limsup_{k\geq \kbar, k \to \infty} \tfrac{\kappa_y}{k-\kbar+1} \sum_{i=\kbar}^k \E[\|X_i-x_\star\|_2],
  \end{align*}
  as desired.
\qed

The conditions in \eqref{eq.noise_avg_ass} in Lemma~\ref{lem.noise_avg} are admittedly nontrivial, but upon closer inspection it is not surprising that conditions akin to these would be needed in order to draw any useful conclusions about an averaged Lagrange multiplier sequence.  First, consider conditions \eqref{eq.yavg_ass_1} and \eqref{eq.yavg_ass_4}, which essentially require that the expected value of the norm of $M_k \Delta_k$ is bounded uniformly over $k \in \N{}$.  Indeed, observe that if one were instead to assume that
\bequation \label{eq.bounded_Mdelta}
  \E[ \|M_k \Delta_k\|_2^2 ] < \infty\ \text{for all}\ k \in \N{},
\eequation
then \eqref{eq.yavg_ass_1} follows directly while one also finds
\begin{align*}
  \sup_{k \in \N{}} \E \left[ \left\|\sum_{i=1}^k \tfrac{1}{\sqrt k} M_i \Delta_i \right\|_2^2 \right] \leq \sup_{k \in \N{}} \tfrac{1}{k} \sum_{i=1}^k \E \left[ \left\|M_i \Delta_i \right\|_2^2 \right] < \infty,
\end{align*}
which yields \eqref{eq.yavg_ass_4} for $\Theta=2$.  It is not surprising that this expected value would be required to be bounded uniformly since the elements of $\{M_k\Delta_k\}$ arise as the error terms in the proof of Theorem~\ref{cor.y_conv_avg}.  Second, beyond merely being bounded, which would not be sufficient on its own, conditions \eqref{eq.yavg_ass_2} and \eqref{eq.yavg_ass_3} require convergence properties of $\{M_k\Delta_k\}$ that match those needed for the multidimensional martingale central limit theorem, which we have stated as Lemma~\ref{lem.mart_clt}.  Notice that if one considers the conditions of Assumption~\ref{ass.M}, which have that the elements of $\{M_k\}$ are defined by a continuous (in fact, Lipschitz continuous) mapping $\Mcal$ of the elements of $\{X_k\}$, then \eqref{eq.yavg_ass_2} and \eqref{eq.yavg_ass_3} essentially require convergence in probability of averages of inner and outer products of the elements of $\{\Mcal(X_k) (\nabla f(X_k) - G_k)\}$, which is a reasonable assumption to consider when, say, $\{X_k\} \asto x_\star$.

We close this section with the following corollary to Theorem~\ref{cor.y_conv_avg}.

\begin{corollary}\label{cor.y_conv}
  Suppose that the conditions of Theorem~\ref{cor.y_conv_avg} hold and that the iterate sequence converges almost surely to $x_\star$, i.e., $\{X_k\} \asto x_\star$.  Then,
  \bequationNN
    \{Y_k^\true\} \asto y_\star\ \ \text{and}\ \ \{Y_k^\avg\}_{k=\kbar}^\infty \asto y_\star.
  \eequationNN
\end{corollary}
\proof
  The almost-sure convergence of $\{X_k\}$ to $x_\star$ implies that
  \bequationNN
    \P \left[ \omega \in \Omega : \limsup_{k \geq \kbar, k \to \infty}\|X_k(\omega) - x_\star\|_2 > 0 \right] = 0.
  \eequationNN
  In addition, from~\eqref{eq.y_conv_true}, one finds for all $\omega \in \Omega$ that
  \bequationNN
    \limsup_{k \geq \kbar, k \to \infty}\|Y_k^\true(\omega) - y_\star\|_2 \leq \limsup_{k \geq \kbar, k \to \infty} \kappa_y \|X_k(\omega) - x_\star\|_2.
  \eequationNN
  These facts combined with nonnegativity of norms yields
  \bequationNN
    \baligned
      &\ \P \left[ \omega \in \Omega : \limsup_{k \geq \kbar, k \to \infty} \|Y_k^\true(\omega) - y_\star\|_2 > 0 \right] \\
      \leq &\ \P \left[ \omega \in \Omega : \limsup_{k \geq \kbar, k \to \infty} \|X_k(\omega) - x_\star\|_2 > 0 \right] = 0,
    \ealigned
  \eequationNN
  which again with nonnegativity of norms further implies
  \bequationNN
    \baligned
      &\ \P \left[ \omega \in \Omega : \limsup_{k \geq \kbar, k \to \infty} \|Y_k^\true(\omega) - y_\star\|_2 = 0 \right] = 1 \\
      \implies &\ \P \left[ \omega \in \Omega : \lim_{k \geq \kbar, k \to \infty} \|Y_k^\true(\omega) - y_\star\|_2 = 0 \right] = 1.
    \ealigned
  \eequationNN
  Consequently, $\{Y_k^\true\}_{k=\kbar}^\infty \asto y_\star$, so $\{Y_k^\true\} \asto y_\star$, as desired.

  Let us now consider the averaged sequence.  Recall that under the conditions of the corollary one has for all $\omega \in \Omega$ that
  \begin{align*}
    &\ \sup_{k \geq \kbar, k \in \N{}} \tfrac{1}{\sqrt{k-\kbar+1}} \left\| \sum_{i=\kbar}^k M_i(\omega) \Delta_i(\omega) \right\|_2 < \infty \\ 
    \implies &\ \limsup_{k \geq \kbar, k \to \infty} \tfrac{1}{k-\kbar+1} \left\| \sum_{i=\kbar}^k M_i(\omega) \Delta_i(\omega) \right\|_2 = 0. 
  \end{align*}
  Thus, one obtains along with \eqref{eq.yavg_ass_4} that
  \begin{align*}
    &\ \P \left[ \omega \in \Omega : \limsup_{k \geq \kbar, k \to \infty} \tfrac{1}{k-\kbar+1} \left\| \sum_{i=\kbar}^k M_i(\omega) \Delta_i(\omega) \right\|_2 = 0 \right] \\ 
    \geq &\ \P \left[ \omega \in \Omega : \sup_{k \geq \kbar, k \in \N{}} \tfrac{1}{\sqrt{k-\kbar+1}} \left\| \sum_{i=\kbar}^k M_i(\omega) \Delta_i(\omega) \right\|_2 < \infty\right] = 1.
  \end{align*}
  Thus, there exists a null (possibly empty) set $\Omega_0 \subset \Omega$ such that
  \bequation\label{eq.y_conv_avg_1}
    \limsup_{k \geq \kbar, k \to \infty} \tfrac{1}{k-\kbar+1} \left\| \sum_{i=\kbar}^k M_i(\omega) \Delta_i(\omega) \right\|_2 = 0 \ \ \text{for all} \ \omega \in \Omega \setminus \Omega_0.
  \eequation
  On the other hand, $\{X_k\} \asto x_\star$ implies $\{\|X_k-x_\star\|_2\} \asto 0$.  Thus, there exists a null (possibly empty) set $\Omega_1 \subset \Omega$ such that
  \bequation\label{eq.y_conv_avg_2}
    \limsup_{k \geq \kbar, k \to \infty} \tfrac{1}{k-\kbar+1} \sum_{i=\kbar}^k \|X_i(\omega) - x_\star\|_2 = 0 \ \ \text{for all} \ \omega \in \Omega \setminus \Omega_1.
  \eequation
  Then, recall from~\eqref{eq.y_conv_avg_pf} that for all $\omega \in \Omega$ one finds
  \begin{align*}
    &\ \limsup_{k \geq \kbar, k \to \infty} \|Y_k^\avg(\omega) - y_\star\|_2 \\ 
    \leq &\ \limsup_{k \geq \kbar, k \to \infty} \left(\tfrac{\kappa_y}{k-\kbar+1} \sum_{i=\kbar}^k \|X_i(\omega) - x_\star\|_2 + \tfrac{1}{k-\kbar+1} \left\|\sum_{i=\kbar}^k M_i(\omega) \Delta_i(\omega) \right\|_2 \right),
  \end{align*}
  which along with the fact that $\Omega_0$ and $\Omega_1$ are null sets, \eqref{eq.y_conv_avg_1}, and~\eqref{eq.y_conv_avg_2} implies
  \begin{align*}
    &\ \P \left[ \omega \in \Omega : \limsup_{k \geq \kbar, k \to \infty} \|Y_k^\avg(\omega) - y_\star\|_2 = 0 \right] \\
    = &\ \P \left[\omega \in \Omega \setminus (\Omega_0 \cup \Omega_1) : \limsup_{k \geq \kbar, k \to \infty} \|Y_k^\avg(\omega) - y_\star\|_2 = 0 \right] \\ 
    \geq &\ \P \Bigg[\omega \in \Omega \setminus (\Omega_0 \cup \Omega_1) : \limsup_{k \geq \kbar, k \to \infty} \Bigg( \tfrac{\kappa_y}{k-\kbar+1} \sum_{i=\kbar}^k \|X_i(\omega) - x_\star\|_2 \\
    &\hspace{150pt} + \tfrac{1}{k-\kbar+1} \Bigg\| \sum_{i=\kbar}^k M_i(\omega) \Delta_i(\omega) \Bigg\|_2 \Bigg) = 0 \Bigg] = 1,
  \end{align*}
  Combined with nonnegativity of norms, it follows that
  \bequationNN
    \P \left[ \omega \in \Omega : \lim_{k \geq \kbar, k \to \infty} \|Y_k^\avg(\omega) - y_\star\|_2 = 0 \right] = 1,
  \eequationNN
  which further implies that $\{Y_k^\avg\}_{k=\kbar}^\infty \asto y_\star$, as desired.
\qed

\section{Numerical demonstrations}\label{sec.numerical}

In this section, we provide the results of a small set of numerical experiments to demonstrate the results of our theoretical analysis.  Specifically, we show when solving each instance in a set of test problems that the primal iterates generated by an implementation of Algorithm~\ref{alg.sqp} approach a primal stationary point (in fact, a global minimizer), and as this occurs one finds that the true Lagrange multiplier sequence converges and averaged Lagrange multiplier sequences converge despite the fact that the sequence of Lagrange multipliers itself does not converge (due to gradient errors).  This behavior of the algorithm is representative of behavior that we have witnessed when solving other problems (not shown due to space considerations) as well.  In terms of averaged Lagrange multipliers, we consider averages over all iterations as well as averages determined through practical strategies that attempt to ignore multipliers obtained at primal iterates that are far from a primal stationary point (or at least the current iterate generated by the algorithm).

For our numerical demonstrations, we consider constrained logistic regression; see also \cite{BCOR23}. In particular, we consider problem instances of the form
\bequationNN
  \min_{x\in\R{n}} \tfrac{1}{N} \sum_{i=1}^N \log (1+e^{-{\gamma_i d_i^Tx}})\ \ \st\ \ Ax = b,\ \ \|x\|_2^2 = 1,
\eequationNN
where $D = [d_1\ \cdots\ d_N] \in \R{n \times N}$ is a feature matrix, $\gamma \in \R{N}$ is a label vector, $A \in \R{m\times n}$, and $b \in \R{m}$.  This problem is nonconvex due to the nonlinear constraint $\|x\|_2^2 = 1$, although we confirmed in all runs that the primal iterates approach the globally optimal solution; see below.  For the feature matrices and label vectors, we considered four datasets from LIBSVM \cite{CL11}, namely, \texttt{a9a} ($(n, N) = (123, 32561)$), \texttt{australian} ($(n, N) = (14, 690)$), \texttt{ionosphere} ($(n, N) = (34, 351)$), and \texttt{splice} ($(n, N) = (60, 1000)$).  These datasets were selected as ones that led to relatively large errors in the stochastic gradient estimates.  (See below for further information about how the gradient estimates were computed in our experiments.)  For each problem instance, we choose $m=10$ and randomly generated the initial point $x_1$ and constraint data $(A,b)$, keeping these quantities fixed for all runs for a given problem instance.  Specifically, each entry of these quantities was drawn from a standard normal distribution and it was confirmed in each case that $A$ had full row rank.

For the implementation of Algorithm~\ref{alg.sqp}, we used the Matlab implementation for \cite{BCRZ21} known as \texttt{StochasticSQP}.\footnote{\url{https://github.com/frankecurtis/StochasticSQP}}  For simplicity, we set $h_k$ as the identity matrix $I$ for all iterations in all runs of the algorithm. We solved each linear system \eqref{eq.system} with Matlab's built-in \texttt{mldivide($\backslash$)} function, which in particular meant that---as would be recommended in general for implementations of the algorithm---all Lagrange multipliers were computed directly from this linear system, not by computing a null space basis $Z_k$.  We chose $\beta_k = \Ocal(\tfrac{1}{k})$ with $\beta_1 = 1$.  For each problem instance, to obtain $(x_\star,y_\star)$, we initially ran \texttt{StochasticSQP} with true gradients in all iterations until a stationary point was found to high accuracy, then verified that at these points both first- and second-order stationary conditions held to high accuracy, thus verifying (approximate) global optimality of the final iterate.  Then, for each instance, we ran \texttt{StochasticSQP} with mini-batch gradient estimates with mini-batch sizes of 16 to generate $\{(x_k,y_k,y_k^\true)\}$.  (For each run, the sequence $\{\beta_k\}$ was chosen as an unsummable, yet square-summable sequence that was tuned in a manner such that the primal iterates tended toward the globally optimal primal solution.)  We verified that in all runs the adaptive parameter rules implemented in \texttt{StochasticSQP} kept the merit parameter sequence $\{\tau_k\}$ constant at the initial value of 0.1 and kept the ratio parameter sequence $\{\xi_k\}$ constant at the initial value of 1.0; hence, the runs respected the simplified variant of the algorithm from \cite{BCRZ21} that we have analyzed as Algorithm~\ref{alg.sqp}.

\begin{figure}
  \centering
  \includegraphics[width=5cm]{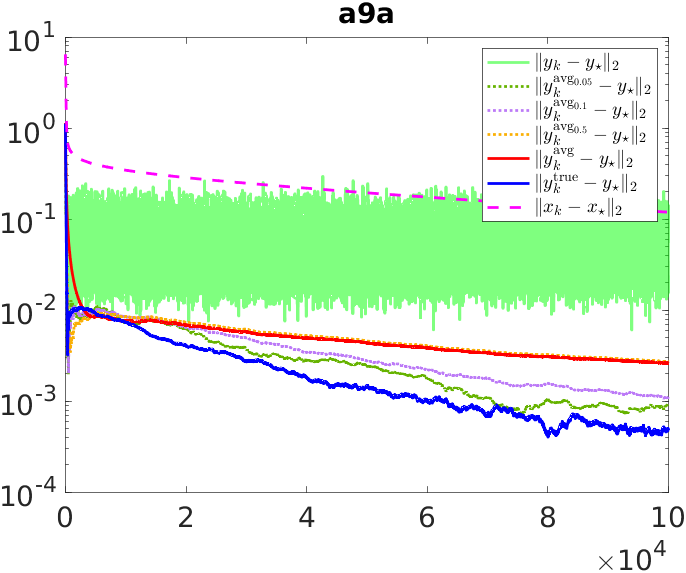}
  \includegraphics[width=5cm]{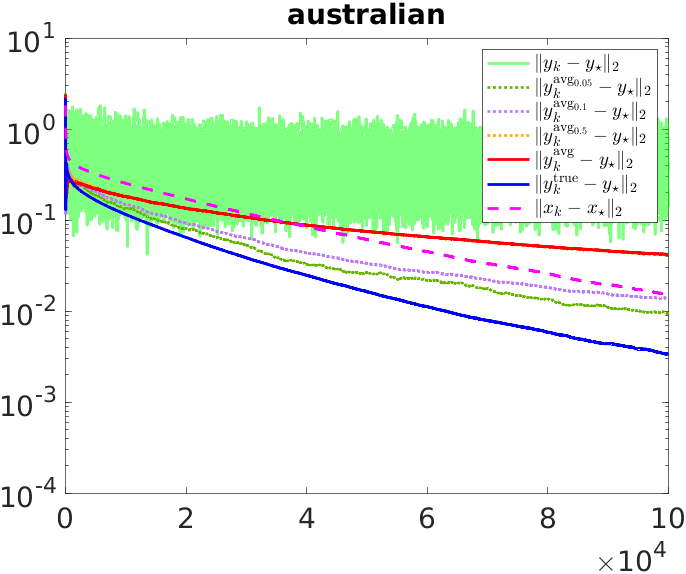}\\
  \includegraphics[width=5cm]{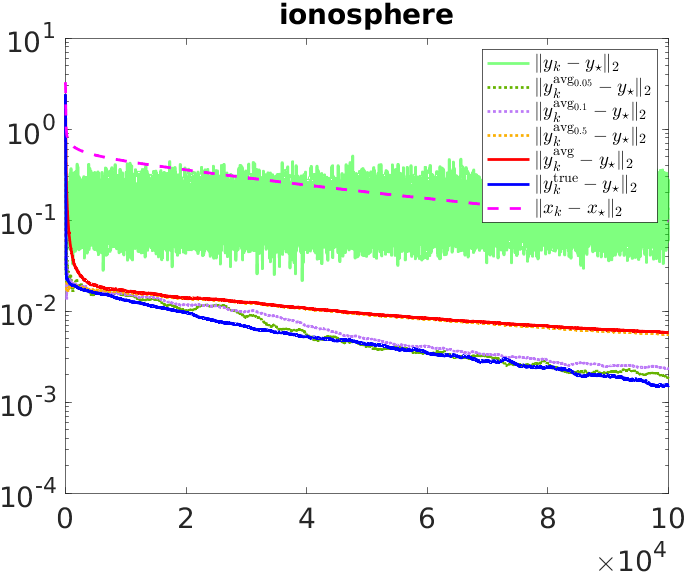}
  \includegraphics[width=5cm]{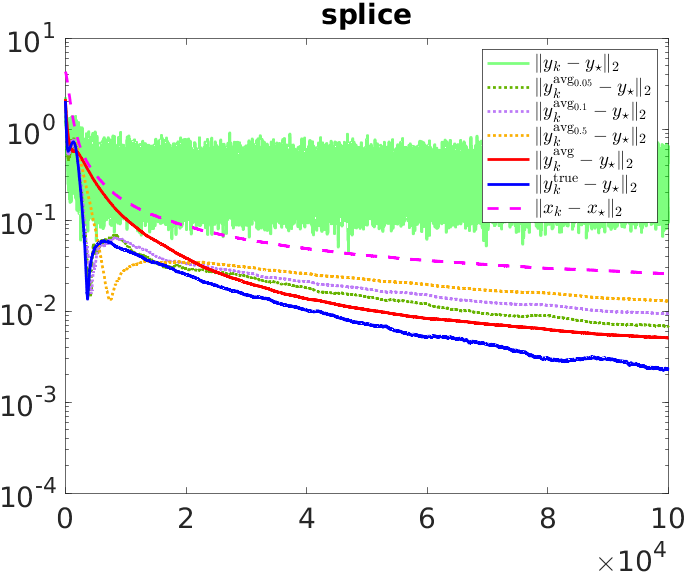}
  \caption{Distances of primal iterates and Lagrange multiplier estimates to solution values over runs of Algorithm~\ref{alg.sqp} to solve constrained logistic regression problems defined using four LIBSVM datasets.}
  \label{fig.dist}
\end{figure}

The results for each problem instance are shown in Figure~\ref{fig.dist}.  Each plot in the figure shows the results of a single run with a budget of $k = 10^5$ iterations, although we verified that the results were qualitatively consistent in other runs (with different random seeds) as well.  In each plot, the figure shows the sequences $\{\|x_k - x_\star\|_2\}$, $\{\|y_k - y_\star\|_2\}$, and $\{\|y_k^\true - y_\star\|_2\}$, where the quantities $x_k$, $y_k$, and $y_k^\true$ are defined as in the prior sections.  The figure also shows the sequences $\{\|y_k^\avg - y_\star\|_2\}$ and $\{\|y_k^{\avg_\epsilon} - y_\star\|_2\}$ for a few values of~$\epsilon \in \R{}_{>0}$, where $y_k^\avg$ is the averaged Lagrange multiplier over all iterations (from the initial point) whereas $y_k^{\avg_\epsilon}$ is the averaged Lagrange multiplier over iterations $\{k',\dots,k\}$, where for each $k \in \N{}$ the index $k'$ is the smallest value in $[k]$ such that $\|x_j - x_k\|_2 \leq \epsilon$ for all $j \in \{x_{k'},\dots,x_k\}$.  In this manner, for relatively large values of $\epsilon$, the average is likely to be taken over more iterations (thus reducing noise), whereas for relatively small values of $\epsilon$, the average is taken only over iterations in which the primal point is relatively close to the current primal iterate (to improve accuracy of the estimate as $x_k$ nears $x_\star$).

Figure \ref{fig.dist} demonstrates the results of our theoretical analysis.  In each run, within the iterations performed before the computational budget was exceeded, the distance of the primal iterate to the primal solution typically reduces monotonically (due to the tuned step-size sequence, as mentioned).  Thus, the run provides an instance in which the primal iterates appear to converge despite the use of stochastic gradient estimates.  At the same time, the distance of the Lagrange multiplier estimate to the optimal Lagrange multiplier oscillates between relatively large values due to the errors in the stochastic gradient estimates, whereas the distances of true and averaged Lagrange multiplier estimates reach much smaller values.  Since the true Lagrange multiplier estimates are unobtainable in practice, one finds that averaged Lagrange multipliers may be considered to obtain better estimates than the $\{y_k\}$ sequence itself.  In some cases, averaging over all iterations may be sufficient, but if desired one may consider choosing $\epsilon \in \R{}_{>0}$ (tuning the value, if computationally feasible) and considering a sequence such as $\{y_k^{\avg_\epsilon}\}$ instead.

\section{Conclusion}\label{sec.conclusion}

We have presented new convergence analyses for a stochastic SQP method built on the stochastic gradient methodology.  In particular, for a simplified variant of the algorithm from \cite{BCRZ21}, we have proved almost-sure convergence guarantees for the primal iterates, Lagrange multiplier estimates, and stationarity measures.  It has been shown that the error in the Lagrange mutipliers is bounded by the distance of the primal iterate to a primal stationary point plus the error due to stochastic gradient estimate.  Furthermore, under modest assumptions, the latter error can be shown to vanish by incorporating a running average of the Lagrange multipliers.  The results of numerical experiments demonstrate our proved theoretical guarantees.

\section*{Acknowledgements}

This work is supported by the U.S.~National Science Foundation (NSF) award IIS-1704458 and U.S.~Office of Naval Research award N00014-21-1-2532.

\bibliographystyle{plain}
\bibliography{ref}

\begin{thebibliography}{10}

\bibitem{BBZ23}
Albert~S. Berahas, Raghu Bollapragada, and Baoyu Zhou.
\newblock An adaptive sampling sequential quadratic programming method for equality constrained stochastic optimization.
\newblock {\em arXiv e-prints}, arXiv:2206.00712v2, 2023.

\bibitem{BCOR23}
Albert~S. Berahas, Frank~E. Curtis, Michael~J. O'Neill, and Daniel~P. Robinson.
\newblock A stochastic sequential quadratic optimization algorithm for nonlinear equality constrained optimization with rank-deficient {Jacobians}.
\newblock {\em arXiv e-prints}, arXiv:2106.13015v2, 2023.

\bibitem{BCRZ21}
Albert~S. Berahas, Frank~E. Curtis, Daniel~P. Robinson, and Baoyu Zhou.
\newblock Sequential quadratic optimization for nonlinear equality constrained stochastic optimization.
\newblock {\em SIAM Journal on Optimization}, 31(2):1352--1379, 2021.

\bibitem{BSYZ23}
Albert~S. Berahas, Jiahao Shi, Zihong Yi, and Baoyu Zhou.
\newblock Accelerating stochastic sequential quadratic programming for equality constrained optimization using predictive variance reduction.
\newblock {\em Computational Optimization and Applications}, 86(1):79--116, 2023.

\bibitem{BXZ23}
Albert~S. Berahas, Miaolan Xie, and Baoyu Zhou.
\newblock A sequential quadratic programming method with high probability complexity bounds for nonlinear equality constrained stochastic optimization.
\newblock {\em arXiv e-prints}, arXiv:2301.00477, 2023.

\bibitem{Ber98}
Dimitri~P. Bertsekas.
\newblock {\em Network Optimization: Continuous and Discrete Methods}.
\newblock Optimization and Neural Computation Series. Athena Scientific, 1998.

\bibitem{BT00}
Dimitri~P. Bertsekas and John~N. Tsitsiklis.
\newblock Gradient convergence in gradient methods with errors.
\newblock {\em SIAM Journal on Optimization}, 10(3):627--642, 2000.

\bibitem{Bet10}
John~T. Betts.
\newblock {\em Practical Methods for Optimal Control and Estimation Using Nonlinear Programming}.
\newblock Advances in Design and Control. Society for Industrial and Applied Mathematics, 2nd edition, 2010.

\bibitem{BCN18}
Léon Bottou, Frank~E. Curtis, and Jorge Nocedal.
\newblock Optimization methods for large-scale machine learning.
\newblock {\em SIAM Review}, 60(2):223--311, 2018.

\bibitem{BN91}
Richard~H. Byrd and Jorge Nocedal.
\newblock An analysis of reduced {Hessian} methods for constrained optimization.
\newblock {\em Mathematical Programming}, 49:285--323, 1990.

\bibitem{BS86}
Richard~H. Byrd and Robert~B. Schnabel.
\newblock Continuity of the null space basis and constrained optimization.
\newblock {\em Mathematical Programming}, 35:32--41, 1986.

\bibitem{CP16}
Antonin Chambolle and Thomas Pock.
\newblock On the ergodic convergence rates of a first-order primal–dual algorithm.
\newblock {\em Mathematical Programming}, 159(1-2):253--287, 2016.

\bibitem{CL11}
Chih-Chung Chang and Chih-Jen Lin.
\newblock {LIBSVM}: {A} library for support vector machines.
\newblock {\em ACM Transactions on Intelligent Systems and Technology}, 2(3):1--27, 2011.

\bibitem{Chu03}
Kai~Lai Chung.
\newblock {\em A Course in Probability Theory, Revised Edition.}, volume 3rd ed.
\newblock Academic Press, 2001.

\bibitem{CDG+22}
Salvatore Cuomo, Vincenzo~Schiano Di~Cola, Fabio Giampaolo, Gianluigi Rozza, Maziar Raissi, and Francesco Piccialli.
\newblock Scientific machine learning through physics-informed neural networks: Where we are and what's next.
\newblock {\em Journal of Scientific Computing}, 92(3):88, 2022.

\bibitem{COR23}
Frank~E. Curtis, Michael~J. O'Neill, and Daniel~P. Robinson.
\newblock Worst-case complexity of an {SQP} method for nonlinear equality constrained stochastic optimization.
\newblock {\em Mathematical Programming}, 2023.
\newblock https://doi.org/10.1007/s10107-023-01981-1.

\bibitem{CRZ23}
Frank~E. Curtis, Daniel~P. Robinson, and Baoyu Zhou.
\newblock Sequential quadratic optimization for stochastic optimization with deterministic nonlinear inequality and equality constraints.
\newblock {\em arXiv e-prints}, arXiv:2302.14790, 2023.

\bibitem{Dur19}
Richard Durrett.
\newblock {\em Probability: Theory and Examples}.
\newblock Cambridge Series in Statistical and Probabilistic Mathematics. Cambridge University Press, 5th edition, 2019.

\bibitem{FNMK22}
Yuchen Fang, Sen Na, Michael~W. Mahoney, and Mladen Kolar.
\newblock Fully stochastic trust-region sequential quadratic programming for equality-constrained optimization problems.
\newblock {\em arXiv e-prints}, arXiv:2211.15943, 2022.

\bibitem{HH14}
Peter Hall and Christopher~C. Heyde.
\newblock {\em Martingale Limit Theory and Its Application}.
\newblock Academic Press, 2014.

\bibitem{JV23}
Xin Jiang and Lieven Vandenberghe.
\newblock Bregman three-operator splitting methods.
\newblock {\em Journal of Optimization Theory and Applications}, 196(3):936--972, 2023.

\bibitem{KKL+21}
George~E. Karniadakis, Ioannis~G. Kevrekidis, Lu~Lu, Paris Perdikaris, Sifan Wang, and Liu Yang.
\newblock Physics-informed machine learning.
\newblock {\em Nature Reviews Physics}, 3(6):422--440, 2021.

\bibitem{KS92}
F.-S. Kupfer and Ekkehard~W. Sachs.
\newblock Numerical solution of a nonlinear parabolic control problem by a reduced {SQP} method.
\newblock {\em Computational Optimization and Applications}, 1(1):113--135, 1992.

\bibitem{LXY23}
Tiejun Li, Tiannan Xiao, and Guoguo Yang.
\newblock Revisiting the central limit theorems for the {SGD}-type methods.
\newblock {\em arXiv e-prints}, arXiv:2207.11755, 2023.

\bibitem{McL74}
Don~L. McLeish.
\newblock Dependent central limit theorems and invariance principles.
\newblock {\em The Annals of Probability}, 2(4), 1974.

\bibitem{NAK23}
Sen Na, Mihai Anitescu, and Mladen Kolar.
\newblock An adaptive stochastic sequential quadratic programming with differentiable exact augmented {Lagrangians}.
\newblock {\em Mathematical Programming}, 199(1-2):721--791, 2023.

\bibitem{NAK23a}
Sen Na, Mihai Anitescu, and Mladen Kolar.
\newblock Inequality constrained stochastic nonlinear optimization via active-set sequential quadratic programming.
\newblock {\em Mathematical Programming}, 2023.
\newblock https://doi.org/10.1007/s10107-023-01935-7.

\bibitem{NM22}
Sen Na and Michael~W. Mahoney.
\newblock Asymptotic convergence rate and statistical inference for stochastic sequential quadratic programming.
\newblock {\em arXiv e-prints}, arXiv:2205.13687, 2022.

\bibitem{Nes05}
Yurii Nesterov.
\newblock Smooth minimization of non-smooth functions.
\newblock {\em Mathematical Programming}, 103(1):127--152, 2005.

\bibitem{OBN21}
Figen Oztoprak, Richard Byrd, and Jorge Nocedal.
\newblock Constrained optimization in the presence of noise.
\newblock {\em arXiv e-prints}, arXiv:2110.04355, 2021.

\bibitem{Poly63}
B.~T. Polyak.
\newblock Gradient methods for minimization of functionals.
\newblock {\em USSR Computational Mathematics and Mathematical Physics}, 3(3):643--653, 1963.

\bibitem{PJ92}
Boris~T. Polyak and Anatoli~B. Juditsky.
\newblock Acceleration of stochastic approximation by averaging.
\newblock {\em SIAM Journal on Control and Optimization}, 30(4):838--855, 1992.

\bibitem{QK23}
Songqiang Qiu and Vyacheslav Kungurtsev.
\newblock A sequential quadratic programming method for optimization with stochastic objective functions, deterministic inequality constraints and robust subproblems.
\newblock {\em arXiv e-prints}, arXiv:2302.07947, 2023.

\bibitem{RDW10}
Tyrone Rees, H.~Sue Dollar, and Andrew~J. Wathen.
\newblock Optimal solvers for {PDE}-constrained optimization.
\newblock {\em SIAM Journal on Scientific Computing}, 32(1):271--298, 2010.

\bibitem{RM51}
Herbert Robbins and Sutton Monro.
\newblock A stochastic approximation method.
\newblock {\em The Annals of Mathematical Statistics}, 22(3):400--407, 1951.

\bibitem{RS71}
Herbert Robbins and David Siegmund.
\newblock A convergence theorem for non-negative almost supermartingales and some applications.
\newblock In {\em Optimizing Methods in Statistics}, pages 233--257. Elsevier, 1971.

\bibitem{Wed73}
Per-Åke Wedin.
\newblock Perturbation theory for pseudo-inverses.
\newblock {\em BIT Numerical Mathematics}, 13(2):217--232, 1973.

\end{thebibliography}

\end{document}